\documentclass[12pt,draftcls,onecolumn]{IEEEtran}
\usepackage[latin1]{inputenc}																							%
\usepackage{type1cm}																										%
\usepackage{type1ec}																										%
\usepackage[T1]{fontenc}																								%
\usepackage{dsfont}																										%
\usepackage{booktabs}																									%
\usepackage{array}																											%
\usepackage{enumerate}																									%
\usepackage[table]{xcolor}																								%
\usepackage{makeidx}																										%
\usepackage{bm}																												%
\usepackage{lastpage}																										%
\usepackage{multicol}																										%
\usepackage{amsmath, amssymb, amsfonts,booktabs}													%
\usepackage{amstext}																										%
\usepackage{ifthen}																											%
\usepackage{ifpdf}																											%
\DeclareMathAlphabet{\mathpzc}{OT1}{pzc}{m}{it}															%
\usepackage{psfrag} 																										%
\usepackage{float} 																											%
\usepackage{cite} 																											%
\usepackage[dvips]{graphicx} 																							%
\usepackage{stackrel} 																										%
\usepackage{graphicx} 																									%
\usepackage{algorithm} 																									%
\usepackage{algorithmic} 																								%
\usepackage{textcomp} 																									%
\usepackage{bbm} 																											%
\usepackage{dsfont} 																										%
\usepackage{subfig} 																										%
\usepackage{multirow}																									%
\usepackage[section,verbose]{placeins} 																			%
\usepackage{tikz}																											%
\usetikzlibrary{arrows,automata,matrix,fit,positioning,calc}											%
\usepackage{pgfplots}																										%
\usepackage[latin1]{inputenc}																							%
\usepackage{verbatim}																									%

\usepackage{color}  
\usepackage[normalem]{ulem} 

\newcommand{\N}{\mathcal{N}}

\usepackage{framed}																							
\usepackage[amsmath,thmmarks,framed,hyperref]{ntheorem}							  %
\newcounter{generalCounter}
\theoremclass		{Theorem}
\theoremstyle		{plain}
\theoremheaderfont	{\normalfont\bfseries}
\theorembodyfont	{\normalfont}
\newtheorem{theorem}{Theorem}
\newtheorem{lemma}[theorem]{Lemma}

\newtheorem{proposition}[theorem]{Proposition}

\usepackage{setspace}

\title{Modular design of jointly optimal controllers and forwarding policies for wireless control }

\title{Modular design of jointly optimal controllers and forwarding policies for wireless control }

\author{\IEEEauthorblockN{Burak Demirel\IEEEauthorrefmark{1}\IEEEauthorrefmark{2}, Zhenhua Zou\IEEEauthorrefmark{1}, Pablo Soldati\IEEEauthorrefmark{3}
and Mikael Johansson\IEEEauthorrefmark{1}}
\thanks{
				\IEEEauthorblockA{\IEEEauthorrefmark{1}ACCESS Linnaeus Centre, KTH Royal Institute of Technology, Osquldas vag~10, SE-100 44 Stockholm, Sweden. }
				\IEEEauthorblockA{\IEEEauthorrefmark{3}Huawei Technologies Sweden AB, Skalhogatan 9-11 box 54, SE-164 94, Kista, Sweden. }
				\IEEEauthorrefmark{2} Corresponding author. E-mail: burak.demirel@ee.kth.se. A preliminary version of parts of this work was presented at several conferences; see~\cite{DZS+:11,Zou12Wiopt,ZouCDC12}.
				}
}
\IEEEoverridecommandlockouts

\begin{document}

\maketitle

\begin{abstract}
We consider the joint design of packet forwarding policies and controllers for wireless control loops where sensor measurements are sent to the controller over an unreliable and energy-constrained multi-hop wireless network. For fixed sampling rate of the sensor, the co-design problem separates into two well-defined and independent subproblems: transmission scheduling for maximizing the deadline-constrained reliability and optimal control under packet loss. We develop optimal and implementable solutions for these subproblems and show that the optimally co-designed system can be efficiently found. Numerical examples highlight the many trade-offs involved and demonstrate the power of our approach.
\end{abstract}


\IEEEpeerreviewmaketitle

\section{Introduction}\label{sec:Introduction}

Cyber-physical systems (CPS) represent a new class of networked embedded systems where requirements on the communication infrastructure are intimately coupled with computational and control requirements. Examples of CPS applications include building management and automation, Internet of things, intelligent transportation systems, and smart grids~\cite{BaG:11,KiK:12}. The industrial practice for designing such systems has been to disregard the computation and communication aspects when designing the controller. This works well if the performance requirements are low, or if there is a significant separation between the time scale of computation and communication, on the one hand, and control, on the other. It is then often possible to make the communication and computation appear reliable and predictable at the time scale of the control loop. When communication, computation and control do interact, the burden of ensuring reliable system-level performance is typically put on the control system. Using high-level abstractions of the deficiencies introduced by unreliable communication and resource-constrained hardware, control algorithms are synthesized to be robust to these uncertainties. However, robustification of control laws can come at a high performance price, and it is sometimes simply not possible to compensate for networking shortcomings in control software.

We argue that efficient CPS systems must be based on the \emph{joint} design of communication, computation, and control. At the same time, it is essential that such a joint design is \emph{modular}, with well-defined interfaces between control algorithms and networking and computation primitives. Modularity allows for specialized development and innovation within each component without affecting the logical correctness of the overall system, and has been a key to massive proliferation in computing and communications. To this end, this paper explores modular co-design of networked control systems with certain optimality properties of the overall system.

Networked control has been an active area of research for more than a decade, and the literature is by now rather extensive, see~\emph{e.g.},~\cite{Nil:98,Hespanha07} and the references therein. The research has mainly focused on control design methods that rely on high-level abstractions of the  communication network in terms of its latency or loss. State-of-the-art control design techniques are extremely powerful when the control system is able to cope with the network deficiencies. However, when the resulting closed-loop performance is unsatisfactory  they typically do not provide any specifics on how the communication system should be modified to yield better performance.  For instance, it is not immediately clear if a shorter sampling interval is better if it also results in a higher packet loss rate in the network. On the communication side, current networking protocols tend to focus on maximizing the long-term throughput or to minimize the energy cost of communication, and do not provide any direct ways to influence control-relevant network performances such as individual packet latency. In fact, little is known about which combinations of end-to-end loss and per-packet deadline guarantees are achievable for a single packet transmitted over an unreliable multi-hop network, and only recently have researchers started to address real-time communication over unreliable networks.

\subsection{Related Work}
Insight into the co-design problem can be obtained from related work in digital control, real-time scheduling, and networked control. We briefly review some of the work that is most relevant for the developments in this paper.

In many cases, faster sampling gives better control performance,  but it also results in controllers that are more ill-conditioned and more sensitive to numerical errors that arise in fixed-point implementations (see~\emph{e.g.},~\cite{FPW:90}). Faster sampling also consumes more computing resources. It is, therefore, often argued that one should choose the longest sampling period that gives acceptable performance. These arguments have~led to a number of well-accepted rules-of-thumb for sample time selection (see~\emph{e.g.},~\cite{AW:96,FPW:90}). In a multi-tasking environment, more periodic tasks could run reliably on the same machine if they run less frequently. The precise number of tasks that can be run depends on what scheduling policy is used. The celebrated schedulability analysis by Liu and Layland~\cite{LiL:73} characterized the maximum utilization for which all tasks can be guaranteed to meet their deadlines under rate-monotonic and earliest-deadline first scheduling. The natural co-design framework is then to adjust the sampling times of controllers to optimize the overall system performance while maintaining schedulability guarantees, see~\emph{e.g.},~\cite{SLSS:96,HeC:05}.

Co-design for networked control systems is more complex, since sampling interval, latency distribution and reliability of end-to-end transmissions all influence the achievable closed-loop performance. Even when these parameters are fixed, the associated optimal control problems have been solved rather recently. Nilsson~\emph{et al.}~\cite{Nil:98} developed linear-quadratic Gaussian optimal control for discrete-time systems with stochastic (networked-induced) delay. Sinopoli~\emph{et al.}~\cite{SSF+:04} focused on the impact of packet losses and developed minimum variance  estimators and characterized their performance, while Imer~\emph{et al.}~\cite{IYB:06} and Schenato~\emph{et al.}~\cite{SSFPS:07} considered linear-quadratic control  under packet losses and established a separation principle under the assumption of reliable and instantaneous acknowledgements from the actuator to controller. The combination of random delays and packet losses was considered by Drew~\emph{et al.} in~\cite{DXG+:05}. Robinson and Kumar~\cite{RoK:08} allowed for uncertain communication also between controller and actuator and studied an optimal controller placement problem (of how to allocate a total uncertainty between the sensor-controller and the controller-actuator communication).

On the networking side, it is well accepted that industrial communication needs reliable packet delivery~\cite{Willig05,Willig08}, but research on wireless communications with hard per-packet deadline constraints is difficult and has appeared only recently. Hou~\emph{et al.}~\cite{HBK:09} proposed a tractable model for deadline-constrained traffic where all packets arrive at
the beginning of an interval and expire at the end of the interval. The main restriction of this work and its extensions~\cite{Hou10a,Hou10b} is that it only deals with a single-hop (star) network topology, while most practical deployments of low-power wireless networks rely on multi-hop communications. In particular, recent communication standards for real-time wireless control, such as \emph{Wireless}HART~\cite{WirelessHARTDataSheet}, ISA-100~\cite{ISA100:09} and IEEE 802.15.4e~\cite{IEEE15.4e:11}, are converging towards a design that combines a multi-hop and multi-path routing with a globally synchronized multi-channel Time Division Multiple Access (TDMA). The communication solutions provided in the present paper target precisely these technologies.
A number of recent papers~\cite{Saifullah10,Munir10} have studied reliable real-time communications over multi-hop networks. Saifullah~\emph{et al.}~\cite{Saifullah10} formulated a scheduling problem for multiple deadline-constrained periodic data flows in \emph{Wireless}HART networks and proved that the problem is NP hard. A heuristic scheduling algorithm without any guarantees on on-time packet delivery was also proposed. A practical scheduling scheme that routes packets on paths with the minimum number of consecutive losses was developed in~\cite{Munir10}.

Several attempts to develop co-design procedures for wireless control systems have appeared in the literature. Early attempts focused on resource-constrained scenarios where the amount of bits that can be communicated over a wireless channel during a sampling interval is limited and needs to be allocated to different control loops~\cite{XJHBG:03}, or assumed that only a single controller can access the communication medium at each sampling instant~\cite{ReS:04}. Liu and Goldsmith~\cite{LiG:04} included detailed models of the communication system, but considered simple network topologies and their designs were neither modular nor optimal. Rabi~\emph{et al.}~\cite{RSPJ:10} focused on co-design of contention-based medium access and networked estimation and studied the interplay between the number of contenders, the sampling interval, and the latency and loss distributions of the sensor-estimator communication. For~\emph{Wireless}HART networks, our earlier paper~\cite{PZSJ:09} argued for structuring the communication schedule into network primitives such as unicast and convergecast, and developed latency-optimal schedules under the assumption that communication links are reliable. In parallel work to this paper, Saifullah~\emph{et al.}~\cite{SWT+:12} developed a heuristic controller-communication co-design approach to calculate sampling intervals of multiple controllers to optimize their overall control performance and ensure schedulability of the real-time communication. The co-design aspect investigated how the additional latency introduced by heuristic retransmission policies (which improve end-to-end reliability) impact the closed-loop performance. Some recent work has studied the energy cost of the wireless network while considering the control performance. Park~\emph{et al.}~\cite{PAJ:11}  tuned protocol parameters to minimize the energy consumption in a wireless network with star-topology while satisfying a desired control performance. Mo~\emph{et al.}~\cite{MGC+:11} later proposed a stochastic sensor scheduling algorithm to minimize the expected estimation error covariance under given energy constraints. In our initial work~\cite{DZS+:11}, we developed a framework to decompose the co-design problem into  well-defined control and communication tasks. We subsequently extended the framework to account also for the energy expenditure in the wireless network~\cite{ZouCDC12}.
\begin{figure}[t]
\centerline{\includegraphics[width=0.6\columnwidth]
  {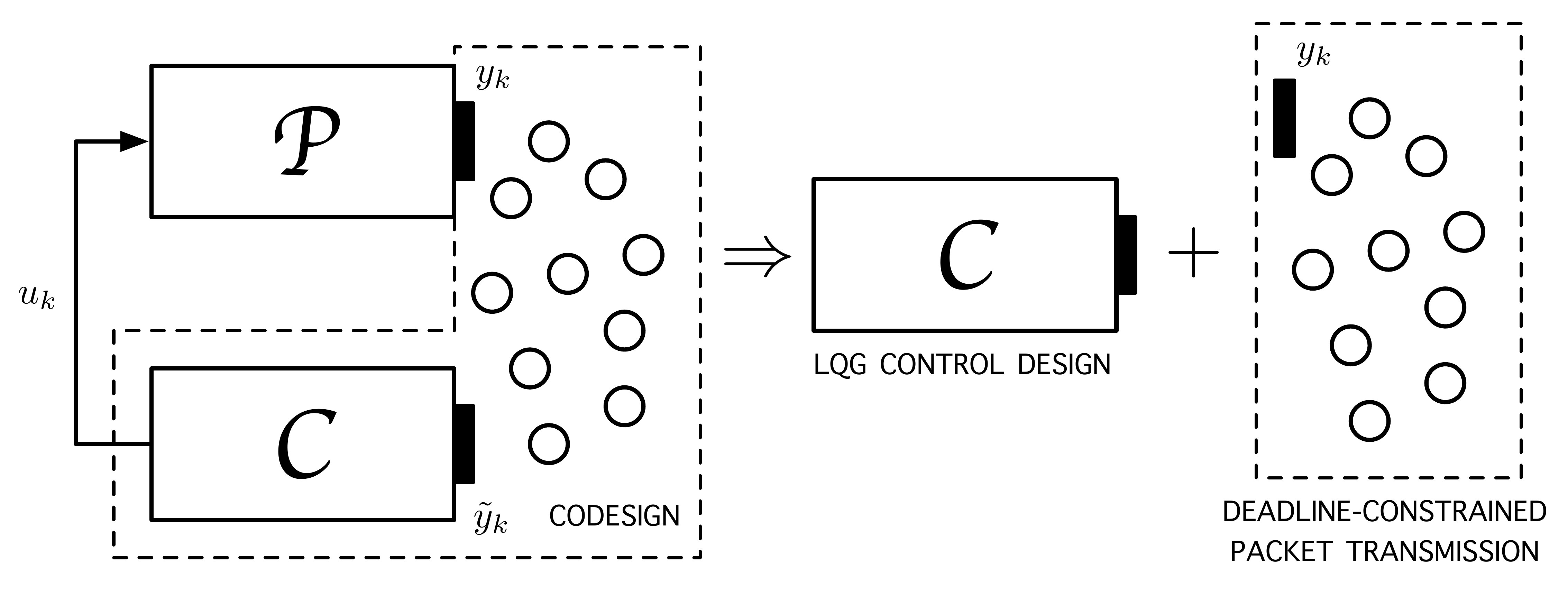}}
  \caption{\normalsize Our co-design framework separates the system design into two well-defined subproblems that admit optimal solutions: deadline-constrained maximum reliability routing, which characterizes the achievable pairs of end-to-end loss and per-packet deadline guarantees,
  and linear-quadratic Gaussian optimal control under latency and loss. The optimal system-level design is obtained by combing these two primitives.}
\label{fig:perfplot}
\end{figure}

\subsection{Contributions}
This paper proposes a co-design framework that finds the jointly optimal controller and multi-hop packet forwarding policy for single-loop wireless control systems. The key is to parameterize the system design in terms of the sample-time of the digital control loop, and to note that the co-design problem then separates into two well-defined sequential design tasks: to schedule the multi-hop network to maximize the deadline-constrained reliability, and to design a controller with optimum performance under (independent) packet losses; see Fig~\ref{fig:perfplot}.
Drawing on our recent work on deadline-constrained scheduling~\cite{SZZJ:10,ZSZJ:10}, we demonstrate how the network scheduling problem can be solved to optimality and how this characterizes the achievable energy-constrained loss-latency region for the multi-hop wireless network. Likewise, for a given communication latency and loss probability, we develop extensions to the work of Schenato~\emph{et al.}~\cite{SSFPS:07} that compute optimal controllers and estimate the associated closed-loop performance. Finally, optimality of the co-design is established by a novel monotonicity result for linear-quadratic control under independent packet loss.
More specifically, the paper contains the following key contributions:
\begin{itemize}
\item We present a co-design framework for wireless control systems where sensor data is forwarded on an unreliable and energy-constrained multi-hop network.
\item By restricting our attention to a time-triggered control architecture, we show that the optimal system performance can be attained by a modular design parameterized by the sampling time of the digital control loop.
\item On the networking side, we derive an optimal multi-hop forwarding policy which maximizes the probability of on-time packet delivery, subject to an energy constraint.
\item On the controller side, we develop the optimal controller under packet losses and characterize its performance. A novel monotonicity result for linear-quadratic control under independent packet losses is established and used to prove optimality of the co-design framework.
\item In numerical examples, we illustrate the power of our framework and explore the trade-offs between sample period, on-time packet delivery probability, network energy consumption and the overall system performance.
\end{itemize}

\subsection{Outline}
The paper is organized as follows. Models and assumptions for the process, sensor, controller, actuator, and network are introduced in~\S\ref{sec:Model_and_Problem_Formulation}.
Then, an optimal and modular co-design framework is proposed in~\S\ref{sec:codesign}. In~\S\ref{sec:lqg_codesign}, the networking and controller subproblems are solved, and optimality of the co-design framework is established. Numerical examples are used to illustrate the power of our framework in~\S\ref{sec:Numerical_Examples}, and conclusions, discussions and directions for future work are stated in~\S\ref{sec:conclusions}. The appendix details proofs of the main theorems.

\subsection{Notation}
In this paper, $\mathbb{N}$ denotes all nonnegative integers, $\mathbb{R}^{n}$ denotes the $n$-dimensional Euclidean space, $\mathbb{R}^{m\times n}$ is the set of all $m\times n$ real
matrices, and $\mathbb{S}_{\geq 0}^{n}~\big(\mathbb{S}_{>0}^{n}\big)$ denotes the cone of real symmetric (positive semi-definite) matrices of dimension $n\times n$. We write the $m\times n$ matrix of all zeros as $0_{m\times n}$, and the $n\times n$ identity matrix as $I_{n\times n}$; $\mathrm{col}\{\lambda_{i}\}$ is the column vector with components $\lambda_{i}$. Let $\textrm{Be}(p)$, $\mathcal{N}(\mu,\sigma^{2})$ and $\textrm{Uni}(a,b)$ denote the Bernoulli distribution, the normal distribution and the uniform (or rectengular) distribution, respectively. Finally, $\mathbbm{1}_{x\in A}$ is the indicator function of the set $A$.


\section{Model and problem formulation}\label{sec:Model_and_Problem_Formulation}

This section summarizes the models and assumptions under which we develop a modular co-design framework with provably optimal performance.

\subsection{Process and sensor} \label{sec:assumptions_sensor}
\begin{figure}\centering
    \includegraphics[angle=0,width=0.7\hsize]{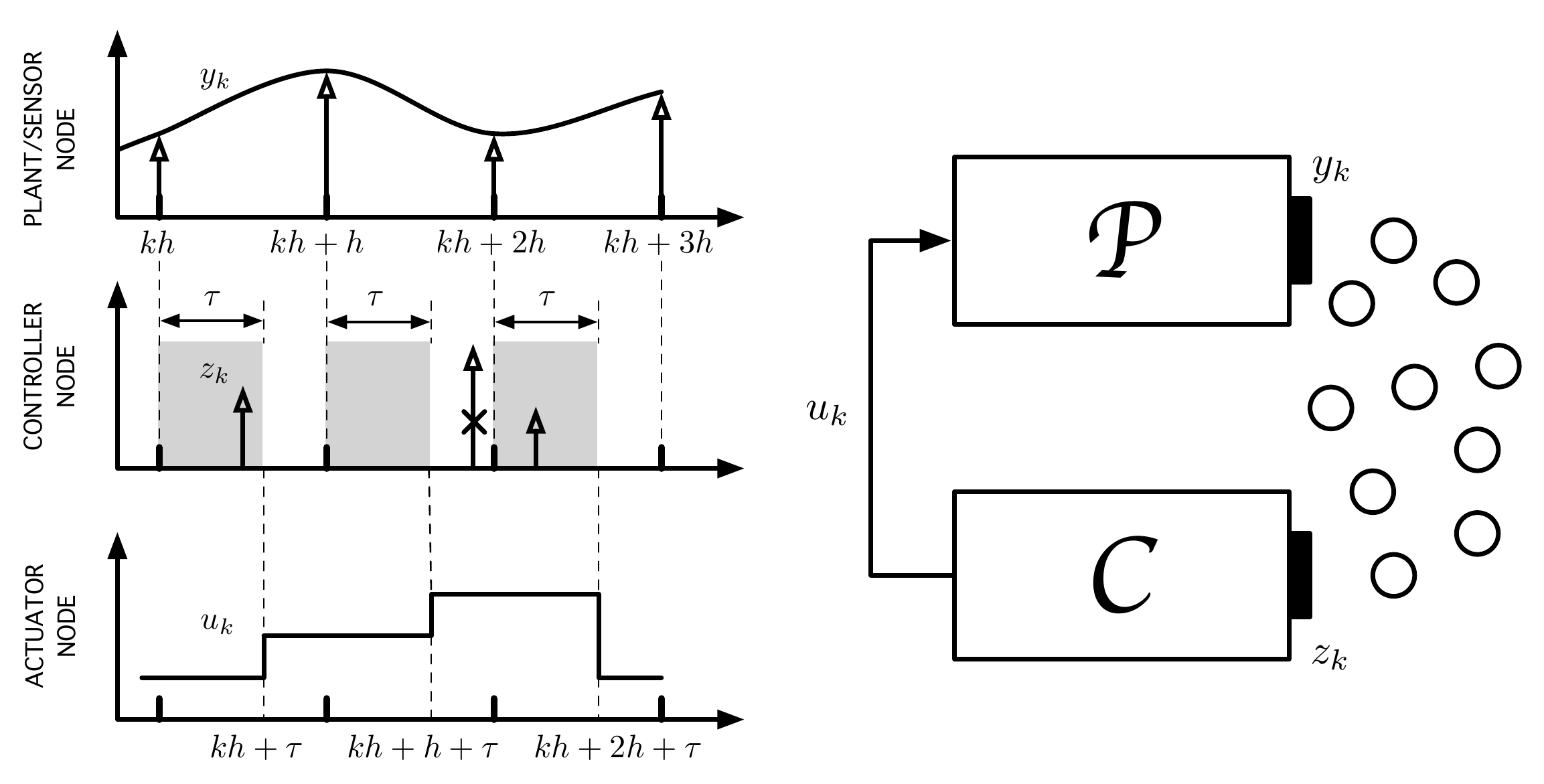}\\
    \caption{Networked control system with timing diagram for sensor, controller and actuator. A sensor takes periodic samples of the system output. The samples are transmitted over a lossy multi-hop wireless network and arrive at the controller with a time-varying delay; packets are dropped if they are unable to meet the per-packet deadline. The controller updates the actuator signal periodically with a fixed lag relative to the time of the sensor readings. }
    \label{fig:TimingDiagram}
\end{figure}

We consider the control of a stochastic linear system
\begin{equation*}
dx(t) = Ax(t)dt + Bu(t)dt + dv_{c}\;, \quad x(0) = x_{0}\;,
\end{equation*}
where $x(t)\in\mathbb{R}^{n}$ is the state, $u(t)\in\mathbb{R}^{m}$ is the control signal, $A\in\mathbb{R}^{n\times n}$ and $B\in\mathbb{R}^{n\times m}$ are the system matrices, and $v_{c}$ is a Wiener process with incremental covariance $R_{v}^{c}\in\mathbb{S}_{\geq 0}^{n}$. Similarly, the initial state $x_{0}$ is modeled as a random variable having a normal distribution with zero mean and covariance $\Sigma_{0}\in\mathbb{S}_{\geq 0}^{n}$, \textit{i.e.}, $x_{0}\thicksim\mathcal{N}(0,\Sigma_{0})$. A noisy measurement of the system output
\begin{equation*}
y(kh) = \tilde{C}x(kh) + w(kh)
\end{equation*}
is taken every sample period $h$. Here, $w(kh)$ is a discrete-time white noise Gaussian process,~independent of the disturbance $v_{c}$, and with zero mean and covariance $R_{w}\in\mathbb{S}_{\geq 0}^{m}$, \textit{i.e.}, $w(kh)\thicksim\mathcal{N}(0,R_{w})$. The sensor measurements are time-stamped and sent over an unreliable multi-hop network as illustrated~in~Fig.~\ref{fig:TimingDiagram}.

\subsection{Controller and actuator} \label{sec:assumptions_control}
We assume the controller and the actuator nodes are co-located, and the control commands are sent from the controller node to the actuator node without information loss or delay. Additionally, the controller and actuator nodes are assumed to be synchronized to the global clock and operate with a fixed lag $\tau \leq h$ relative to the sampling times of the sensor. As shown in Fig.~2, the applied actuator command in each sampling interval is then
\begin{align}
u(t) = \Bigg\{
\renewcommand{\arraystretch}{1.25}
\setlength{\arraycolsep}{1.6pt}
\begin{array}{lrcccl}
u(kh-h), & kh & \leq & t & < & kh+\tau \;, \\
u(kh), & kh+\tau & \leq & t & < & kh+h \;.
\end{array}
\end{align}
In summary, the controller uses the information sensed at time $kh$ and available at time $kh+\tau$ to compute the control action between the controller updates $\big[kh+\tau, (k+1)h+\tau\big)$ for $k \in \mathbb{N}$.

\subsection{Multi-hop wireless network}
The multi-hop wireless network consists of a set of nodes $\mathcal{N}=\{1,\dotsc,Z\}$ equipped with half-duplex radio transceivers. The destination node (collocated with the controller) is labeled $Z$. We represent the network topology as a directed graph $\mathcal{G} = (\mathcal{N},\mathcal{L})$ with nodes ${\mathcal N}$ and links ${\mathcal L}$. The presence of a directed link $(i,j) \in \mathcal{L}$ means that node $i$ is able to deliver a packet successfully to node $j$. Nodes are synchronized to the global clock and communication is slotted. Each slot is $t_s$ milliseconds long and allows for the transmission of a single packet and the reception of the associated  acknowledgement from the next-hop node. Packet transmissions are unreliable and considered successful only if both the packet and the acknowledgement are delivered successfully. We assume independent erasure events following a Bernoulli process with loss probabilities ${\mathbf p}=[p_{ij}]$ (\emph{i.e.}, packet transmission on link $(i,j)$ fails with probability $p_{ij}$, independently of other links). Moreover, nodes do not have access to the current channel state, only their statistics (\textit{i.e.}, loss probabilities on their outgoing links).

The main source of energy consumption in a sensor network is the radio transmission. In this paper, we assume a constraint $\varepsilon$ on the average number of transmission attempts on the wireless sensor network per millisecond in order to guarantee a desired life time. For ease of presentation, we normalize the energy cost of one transmission to one unit of energy.

Based on the time-triggered sensing and control model in \S~\ref{sec:assumptions_sensor} and \S~\ref{sec:assumptions_control}, one packet is injected in the network every sampling period $h$ and must reach the destination node $Z$ within a hard deadline $\tau$, after which it is declared as a loss. The network problem is then to develop a packet forwarding policy $\pi$ that determines if a node should forward a packet or drop it, and to which node it should attempt to transmit, while considering the hard packet deadline and the energy constraint limiting the expected number of transmission attempts.

\subsection{System-level performance and co-design objective}
We aim at developing a multi-hop packet forwarding policy $\pi$, satisfying the energy constraint, and a controller that together minimize the continuous-time closed-loop loss function
\begin{align}
\renewcommand{\arraystretch}{1.0}
\setlength{\arraycolsep}{2.00pt}
J_{c} = \mathbf{E}\bigg\{\int_{0}^{\mathrm{T}}\left[%
\begin{array}{c}
  x(t) \\
  u(t) \\
\end{array}%
\right]^{\intercal}
\left[%
\begin{array}{cc}
  Q_{xx}^{c}  & Q^{c}_{xu} \\
  Q_{xu}^{c\intercal} & Q^{c}_{uu} \\
\end{array}%
\right]
\left[%
\begin{array}{c}
  x(t) \\
  u(t) \\
\end{array}%
\right]dt + x^{\intercal}(\mathrm{T})Q^{c}_{0}x(\mathrm{T})\bigg\}\;, \label{eqn:continuous_loss}
\end{align}
subject to the stochastic system constraints
\begin{align*}
dx = &\; Axdt + Budt + dv_{c}, \\
y(s) = &\; \rho_{k}^{\pi}\tilde{C}x(kh) + w(kh), \;\, kh+\tau\leq s \leq kh+h+\tau
\end{align*}
for $k\in\mathbb{N}$. Here, $Q^{c}_{xx}$ and $Q_{0}^{c}$ are symmetric and positive semi-definite matrices while $Q^{c}_{uu}$ is symmetric and positive definite. The controller computes the control sequence $\{u(k)\}_{k\geq 0}$ based on the available sensor information $y(s)$. Note that sensor data is delayed by a fixed time $\tau$ and may be lost; $\rho^{\pi}_{k}\in \{0,1\}$ is an indicator variable representing whether or not the forwarding policy $\pi$ was able to deliver the $k$-th sensor data packet within the time limits, specified by the hard deadline $\tau$. The expectation in Eq.~\eqref{eqn:continuous_loss} is taken over the random process and measurement disturbances, the initial condition, as well as over the random packet losses on the individual links in the network. We denote the aforementioned optimization problem $\mathbf{OP1}$.

\section{A Modular Co-Design Framework}\label{sec:codesign}

In this section, we will describe the key ideas behind our co-design framework for linear quadratic control over unreliable and energy-constrained multi-hop networks. The framework is modular as it separates the co-design problem into two well-defined design tasks: one for the controller and another one for the network. It is optimal since it allows to find the jointly optimal networking and control design subject to our assumptions and restrictions.

In its general form, the co-design problem is to jointly optimize the closed-loop control loss over the fixed time lag $\tau$ (\emph{i.e.}, the deadline of the packet), the sampling interval $h$ (\emph{i.e.}, the packet generation rate), the packet forwarding policy $\pi$, and the control law $u$ for computing the actuator command. It is important to notice that closed-loop loss is stated in continuous time to allow for comparison of solutions that use different sampling intervals. As we will show in \S~\ref{subsec:control_design}, the fact that the actuator holds the control signals over the intervals $\big[kh+\tau, kh+h\big)$ allows the optimization problem $\mathbf{OP1}$ to be converted to an equivalent discrete-time problem %
\begin{align*}
	\renewcommand{\arraystretch}{1.2}
	\setlength{\arraycolsep}{2.25pt}
  	\begin{array}[c]{ll}
  		\underset{\tau, h, \pi, u}{\mbox{minimize}}
    & J_{d}(h,\tau) \\
  		\mbox{subject to} & \Bigg\{
  		\begin{array}{rll}
			\xi_{k+1} & = & \Phi(h,\tau)\xi_{k} + \Gamma(h,\tau)u_{k} + v_{k} \\
			y_{k} & = & \rho_{k}^{\pi}\tilde{C} \xi_{k} + w_{k}.
		\end{array}
  	\end{array}
\end{align*}

To establish modularity and optimality, we will demonstrate in \S~\ref{sec:optimality} that for a fixed sampling period $h$ and time lag $\tau\leq h$, the optimal control loss is monotone decreasing in the \emph{deadline-constrained reliability} $\rho^{\pi}=\mathbf{E}\{ \rho^{\pi}_k\}$, \emph{i.e.}, the probability that a packet arrives at the controller node within a deadline $\tau$. This implies that for a given $\tau$ and $h$ and under the restriction of one-sample delayed time-triggered control architectures, the optimal co-designed system is obtained by the following two distinguished sequential design problems:
\begin{itemize}
\item[(a)] developing a packet forwarding policy $\pi^\star$ of the network that maximizes the deadline-constrained reliability with deadline of
    $D=\lfloor \frac{\tau}{t_s} \rfloor$ time slots subject to energy constraint,~\textit{i.e.},
\begin{align*}
  \begin{array}[c]{lll}
  \pi^{\star}= & \underset{\pi}{\mbox{argmax}} & \rho^{\pi}\\
  &\mbox{subject to} & C^{\pi} \leq C_{\text{req}},
  \end{array}
\end{align*}
where $C^\pi$ is the average number of transmissions per injected packet attempted by the forwarding policy $\pi$ and $C_{\text{req}}=\varepsilon h$ is the constraint on the average number of transmissions per injected packet due to the energy constraint $\varepsilon$.
\item[(b)] computing the control action using the optimal linear-quadratic controller under packet loss, \textit{i.e.},
\begin{align*}
	\renewcommand{\arraystretch}{1.2}
	\setlength{\arraycolsep}{2.25pt}
  	\begin{array}[c]{ll}
  		\underset{u}{\mbox{minimize}} & J_{d}(h,\tau) \\
  		\mbox{subject to} & \Bigg\{
  		\begin{array}{rll}
			\xi_{k+1} & = & \Phi(h,\tau)\xi_{k}
                        + \Gamma(h,\tau)u_{k} + v_{k} \\
			y_{k} & = & \rho_{k}^{\pi^{\star}} \tilde{C} \xi_{k} + w_{k}
		\end{array}
  	\end{array}
\end{align*}
where $\rho_{k}^{\pi^{\star}}\in \{0,1\}$ is an indicator variable for on-time packet delivery by $\pi^{\star}$. We will show that under Bernoulli link losses and the optimal forwarding policy, we have  $\rho_{k}^{\pi^{\star}}\thicksim \text{Be} (\rho^{\pi^\star})$.
\end{itemize}
Lastly, the optimal solution is obtained by sweeping over all admissible $\tau$ and $h$ values, and identifying the pair of $\tau$ and $h$ values with minimum control loss.

\section{Co-design for linear-quadratic control} \label{sec:lqg_codesign}
We now turn our attention to the development of the optimal control design and forwarding policies required by the co-design framework for linear-quadratic optimal control over a network where packet losses on individual links follow independent Bernoulli processes. The \emph{deadline-constrained maximum-reliability forwarding problem}, defined and solved in \S~\ref{subsec:IVA}, provides the optimal network operation for a fixed sampling interval $h$ and time lag $\tau$. The \emph{optimal control under independent packet losses} is developed in
\S~\ref{subsec:control_design}. Finally, \S~\ref{sec:optimality} establishes~a monotonicity property of the optimal control loss that allows to conclude optimality of the design.

\subsection{Deadline-constrained maximum reliability forwarding} \label{subsec:IVA}

The data flow from sensor to controller node is periodic, with a new packet being produced every $h$ milliseconds. Since the controller disregards packets that have not been delivered within the deadline $\tau \leq h$, nodes can drop packets that will not be able to meet their deadline, which implies that there will be at most one sensor packet in the network at any point in time.  Moreover, since the packet loss processes on links are assumed to be memoryless, the optimal policy for the periodic flow can be constructed by the repeated application of the optimal forwarding policy for a single transient packet.
In the \emph{deadline-constrained} forwarding problem, we pose a strict deadline of $D=\lfloor \frac{\tau}{t_s} \rfloor$ time slots and look for forwarding policies that maximize the probability that the sensor packet is delivered to the destination node (directly connected to the controller) within the deadline. Since the packet is injected every $h$ milliseconds and there is at most one packet on the network during this interval, the energy constraint then can be transformed into the constraint on the average number of transmissions per injected packet $C_{\text{req}} = \varepsilon h$.

More specifically, we consider a scenario where a single packet, generated by the source node at time $t=0$, should be transmitted over a multi-hop wireless network to the destination node
$Z$ within a deadline of $D$ time slots, and the expected number of transmissions of each packet should be smaller than $C_{\text{req}}$. The goal is to derive an optimal forwarding policy $\pi^\star$ with the maximum packet delivery reliability.

\subsubsection{Constrained Markov Decision Process Formulation}
The deadline-constrained packet forwarding problem can be formulated as a finite-horizon Markov decision process (MDP)~\cite{Put:94} with the  horizon equal to the packet deadline $D$.

The state $s(t)$ of the MDP is the packet location at time $t$. The state at the next time $s(t+1)$ depends on the action $a(t)$ that chooses whether or not to attempt a transmission, and if a transmission is attempted, to which next-hop node the packet should be addressed. The state transition probability $\Pr\{s(t+1)|s(t),a(t)\}$, which describes the probability that the next state is $s(t+1)$ given that the current state is $s(t)$ and the current action is $a(t)$, is determined by the loss probability of the associated link. If the action is to hold the packet (\textit{i.e.}, $a(t) = s(t)$), then $s(t+1)=s(t)$. On the other hand, if the action is to forward the packet, then the packet can be either at node $s(t)$ or at neighbor $a(t)$ to which the packet is forwarded. Recall that $p_{ij}$ denotes the loss probability of link $(i,j)$. In summary, the state transition probability is
\begin{align*}
& \Pr \{ s(t+1) | s(t), a(t) \} 
= \begin{cases}
1 & \mbox{if } a(t)=s(t), s(t+1)=s(t), \\
p_{s(t)a(t)} \; & \mbox{if } a(t) \neq s(t), s(t+1)=s(t), \\
1 - p_{s(t)a(t)} \; & \mbox{if } a(t) \neq s(t), s(t+1)=a(t), \\
0 & \mbox{otherwise.}
\end{cases}
\end{align*}

There is no reward for each action, but a terminal reward $\mu_D\bigl(s(D)\bigr)$ is given if the packet is at the destination node $Z$ at the last time slot $D$, and the energy cost $c\bigl(s(t),a(t)\bigr)$ is incurred when the action is to transmit the packet. They are defined separately as
\begin{align*}
\mu_D \bigl(s(D)\bigr) &=
\begin{cases}
1 & \quad \text{if } s(D) =
Z
, \\
0 & \quad \text{otherwise};
\end{cases}
\\
c\bigl(s(t),a(t)\bigr) &=
\begin{cases}
1 & \quad \text{if } a(t) \neq s(t), \\
0 & \quad \text{otherwise}.
\end{cases}
\end{align*}
Note that in the MDP formulation, we enforce the packet to stay at the destination if it arrives before the deadline, and the reward is collected only at the deadline time $t=D$.

The core problem of an MDP is to find a sequence of actions such that the expectation of the sum of the rewards is maximized or the sum of the costs is minimized. We define a decision rule that prescribes the action in each state at a specified time, and a policy composed by a sequence of decision rules at each time. In the most general case, the decision rule should depend on the previous states and actions, and can prescribe actions randomly.
More precisely, let $H(t)$ be the set of all possible histories where a history is a sequence of previous states and actions, \textit{i.e.}, $\bigl(s(0), a(0), \dotsc, s(t-1), a(t-1), s(t)\bigr)$.
The decision rule is a function $d(t) : H(t) \rightarrow \mathcal{P}\bigl(A(t)\bigr)$ that maps $H(t)$ into a set of probability distributions on the action space $A(t)$ of all possible actions. Since the MDP is a process of finite-horizon with length $D$, the policy is $\pi \triangleq \bigl(d(0),d(1), \dotsc, d(D-1)\bigr)$, indexed by time.
Under a policy $\pi$, the expected total reward (deadline-constrained packet
reliability) is defined as
\begin{align*}
\rho^\pi \triangleq \mathbf{E}^\pi_{s(0)}
  \bigl\{ \mu_D \bigl(s(D)\bigr) \bigr\},
\end{align*}
where $s(0)$ is the initial packet location (\textit{i.e.}, the source node).
Note the expectation $\mathbf{E}$ is taken over the probability space induced by policy $\pi$.
Similarly, the expected energy cost of an end-to-end packet delivery on the network is
\begin{align*}
C^\pi \triangleq \mathbf{E}^\pi_{s(0)}
\biggl\{\sum_{t=0}^{t=D-1} c\bigl(s(t),a(t)\bigr) \biggr\}.
\end{align*}

The maximum deadline-constrained reliability with energy cost constraint problem is then
\begin{align}
  \begin{array}[c]{ll}
  \underset{\pi}{\mbox{maximize}} & \rho^{\pi}\\
  \mbox{subject to} & C^{\pi} \leq C_\text{req}.
  \end{array}
  \label{eq:cmdp}
\end{align}

\subsubsection{Randomized policy}
This problem falls into the category of constrained MDP (CMDP). The standard solution if we are only interested in the total reward is to use linear programming~\cite{Put:94}. However, we are also interested in the structure of the optimal policy to understand which forwarding logic to implement in individual nodes. To this end, we use the Lagrangian approach proposed in~\cite{Alt:99} for CMDP to convert it to a non-constrained weighted sum problem.

The Lagrange dual problem of~\eqref{eq:cmdp} is
\begin{align}&
\begin{array}[c]{ll}
\underset{\delta}{\mbox{minimize}} & \quad \max\limits_\pi
\{ \rho^\pi - \delta \cdot C^\pi \} + \delta \cdot C_\text{req}\\
\text{subject to} & \quad \delta \geq 0 \;.
\end{array}
\label{eq:Dual}
\end{align}
Our finite-horizon CMDP can be cast as an infinite-horizon CMDP with total cost criterion. The Markov state $s$ is extended to include the time from $t=0$ to $t=D$ (\emph{i.e.} states are now node-time pairs). It goes to the next state with time $t+1$ only if the current state's time is $t$. We define a termination state to which all the states with time $D+1$ are directed. This is an absorbing state with no reward and cost. All other parameters including rewards, costs and state transition probabilities remain the same. It can be shown that this is a contracting MDP as defined in~\cite[Def.~2.4]{Alt:99}. Hence, by~\cite[Thm.~4.8 ii]{Alt:99}, the duality gap is zero. To solve problem~\eqref{eq:Dual}, we hence need to solve the weighted sum maximization of reliability and energy, \textit{i.e.},
\begin{align}
\max_{\pi} \{ \rho^\pi - \delta \cdot C^\pi \}
\label{eq:WeightedSum}
\end{align}
$\text{for a given } \delta \geq 0$.

By treating the weighted energy cost $\delta \cdot C^\pi$ as a negative reward scaled by $\delta$ in the MDP formulation, a history-independent and deterministic optimal policy for the weighted sum maximization can be found by dynamic programming~\cite{Put:94}.
The details of the dynamic programming framework will be developed in \S~\ref{sec:DPSolution}.
Note that the maximum reliability forwarding problem with no energy constraint is a special case of problem~\eqref{eq:WeightedSum} with $\delta=0$.

We will first show that an optimal policy for the minimum energy problem can be constructed by randomizing between two deterministic policies, each of which is optimal for a different value of $\delta$ in the weighted sum problem.
We further explicitly derive these two policies and specify the probabilities at which they are selected in the optimal randomized policy.

Let $C^\star (\delta)$, $\rho^\star (\delta)$ and $\pi^\star (\delta)$ be the optimal energy, reliability and policy in the weighted sum problem for a given $\delta$, respectively. Define $\mathcal{C} \triangleq \{ C^\star (\delta), \; \text{for all } \delta\}$ and $\Delta_C \triangleq \{ \delta : C^\star(\delta) = C\}$ for a given $C \in \mathcal{C}$. We have the following results:
\begin{lemma}
\label{lem:ParetoFrontier}
$\mathcal{C}$ is a finite set. For a given $C \in \mathcal{C}$, $\rho^\star(\delta)$ is unique for all $\delta \in \Delta_C$.
\end{lemma}
{\noindent \textit{Proof:}}
See Appendix~\ref{sec:proofParetoFrontier}.
\hfill $\blacksquare$

\begin{theorem}
\label{thm:ParetoFrontier}
Let $C^{(1)} = \max \{ C \in \mathcal{C} : C \leq C_\text{req}\}$ and $C^{(2)} = \min \{ C \in \mathcal{C} : C > C_\text{req}\}$ with the associated unique reliability $\rho^{(1)}$ and $\rho^{(2)}$. The maximum reliability of the forwarding problem~is
\begin{align}
\rho^\star = \rho^{(1)} + \dfrac{C_\text{req}-C^{(1)}}{C^{(2)}-C^{(1)}}
 (\rho^{(2)} - \rho^{(1)}).
\label{eq:minEnergySol}
\end{align}
Suppose that the
history-independent and deterministic
optimal policies that attain $(\rho^{(1)}, C^{(1)})$ and $(\rho^{(2)}, C^{(2)})$ are $\pi^{(1)}$ and $\pi^{(2)}$ respectively. An optimal policy $\pi^\star$ for the minimum energy problem is obtained by random selection of policies $\pi^{(1)}$ and $\pi^{(2)}$ with probabilities
\begin{align*}
\theta^{(1)} = \dfrac{C^{(2)}-C_\text{req}}{C^{(2)}-C^{(1)}} \; ; \quad
\theta^{(2)} = \dfrac{C_\text{req}-C^{(1)}}{C^{(2)}-C^{(1)}}.
\end{align*}
\end{theorem}
{\noindent \textit{Proof:}}
Lemma~\ref{lem:ParetoFrontier} shows the existence of $C^{(1)}$ and  $C^{(2)}$ and the uniqueness of $\rho^{(1)}$ and $\rho^{(2)}$.
The rest of the proof is in Appendix~\ref{sec:proofParetoFrontier}.
\hfill $\blacksquare$

The theorem states that the optimal forwarding policy is to make a random selection between two history-independent and deterministic policies, each found by dynamic programming.
Combined with the assumption that link losses are independent, we can conclude that the packet loss by the optimal forwarding policy is a Bernoulli random process.


\subsubsection{General dynamic programming framework}
\label{sec:DPSolution}
Next, we will develop a general dynamic programming framework to solve the weighted sum maximization problem~\eqref{eq:WeightedSum}.

Since the MDP state is composed of the packet location, let the maximum utility at time $t$ and node $i$ be
\begin{align*}
U^\star_i (t) &=\rho^\star_i(t) - \delta C^\star_i(t).
\end{align*}
This quantity describes the optimal utility for packet delivery within the next $D-t$ time slots; $\rho^\star_i(t)$ and $C^\star_i(t)$ are the corresponding optimal reliability and optimal energy cost, respectively.

\begin{figure}[th!]
  \centerline{\includegraphics[width=0.45\hsize]{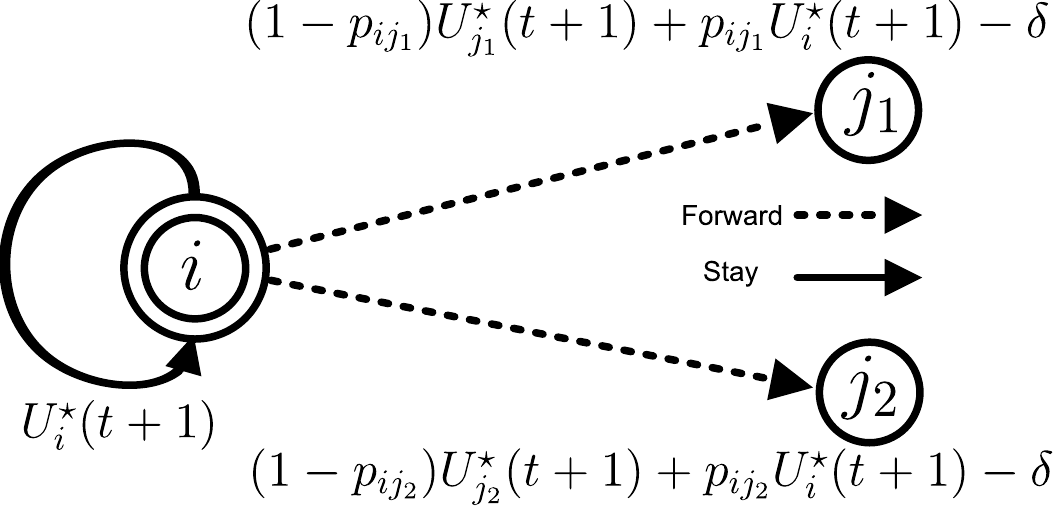}}
  \caption{An illustration of Bellman equation at node $i$ and time $t$ for two outgoing links.}
  \label{fig:forwardEx}
\end{figure}

To this end, our aim is to develop optimal forwarding policies that attain $U_i^{\star}(0)$, the maximum utility of node $i$ at $t=0$ with a deadline of $D$ time slots. This quantity can be computed recursively by dynamic programming from $t=D-1$ to $t=0$, starting from initial condition
\begin{align*}
U^{\star}_i (D) &=
\begin{cases}
1 &  \mbox{if $i =
Z
$}, \\
0 &  \mbox{if $i \neq
Z
$};
\end{cases}
\\
\rho^{\star}_i (D) &=
\begin{cases}
1 &  \mbox{if $i =
Z
$}, \\
0 &  \mbox{if $i \neq
Z
$};
\end{cases}
\quad
C^{\star}_i (D) =0.
\end{align*}
At each step $t < D$, the maximum utility $U^{\star}_i(t)$ at node $i$ is characterized by the Bellman~equation,
\begin{align}
U^{\star}_i (t) = \max\Big\{ \max_{j \in \N_i} U_i^{j} (t) ,\; U_i^i (t) \Big\} \label{eq:MaximumUtility}
\end{align}
where $U_i^{j} (t)$ is the utility of forwarding to a neighbor $j\in\N_i$, and $U_i^i (t)$ is the utility of withholding the packet at node $i$, respectively. These utilities are computed as
\begin{align*}
  U_i^{j} (t) & =
  \underbrace{(1-p_{ij}) U^{\star}_j(t+1)}
   _{\mbox{\small{Success forwarding}}}
  + \underbrace{ p_{ij}  U^{\star}_i(t+1)}
   _{\mbox{\small{Failed forwarding}}}
   - \underbrace{\delta}_{\mbox{\small{Energy cost}}};
  \\
  U_i^i (t) &= U^{\star}_i(t+1).
\end{align*}%
An illustration of the Bellman equation update for two outgoing links are shown in Fig.~\ref{fig:forwardEx}.
Note that in each step $t$, the update in Eq.~\eqref{eq:MaximumUtility} requires only the maximum utility $U^{\star}_i(t+1)$ of node $i$ and the maximum utility $U^{\star}_j(t+1)$ of node $j$ where $j$ is its one-hop away neighbor.

The optimal action at time $t$ forwards the packet to the node $j^{\star}_i(t)$ that obtains the maximum in Eq.~\eqref{eq:MaximumUtility}, and $j^\star_i(t)$ is given by
\begin{align}
\label{eq:OptimalPolicy}
j^{\star}_i(t)=
\begin{cases}
i &
\text{if~} U_i^i (t)
\geq U_i^{j} (t)
\; \forall j \in \mathcal{N}_i; \\
\arg \max\limits_{j \in \mathcal{N}_i}
U_i^{j} (t) &
\text{otherwise.}
\end{cases}
\end{align}%

Note that withholding the packet does not consume energy, and that we break ties arbitrarily.
With Eq.~\eqref{eq:OptimalPolicy}, one can compute $\rho^{\star}_i(t)$ and $C^{\star}_i(t)$. If $j^{\star}_i(t) \neq i$, then
\begin{align*}
\rho^{\star}_i (t)
&= \big(1- p_{ij_i^\star(t)} \big) \rho^{\star}_{j^{\star}_i(t)}(t+1)
 + p_{ij_i^\star(t)} \rho^{\star}_i(t+1); \\
C^{\star}_i (t) &=
\big( 1- p_{ij_i^\star(t)} \big) C^{\star}_{j^{\star}_i(t)}(t+1)
+ p_{ij_i^\star(t)} C^{\star}_i(t+1) + 1.
\end{align*}
On the other hand, if $j^{\star}_i(t)=i$, then
\begin{align*}
\rho^{\star}_i (t) = \rho^{\star}_i(t+1)
\; \quad
C^{\star}_i (t) = C^{\star}_i(t+1).
\end{align*}

The forwarding policy from the dynamic programming framework is composed of the actions at each possible packet location and time. We can implement this policy in a distributed fashion at each node without need for inter-node coordination. A node stores only the policy associated with itself as a lookup table indexed by time. At run time, if a node holds the packet, then it forwards the packet according to the lookup table. We also showed in Theorem~\ref{thm:ParetoFrontier} that the optimal forwarding policy is to make a random selection between two policies found by dynamic programming framework. A naive implementation of the optimal forwarding policy would be to randomly select one of the deterministic policies when the packet is created, mark the packet accordingly, and let intermediate nodes forward according to the chosen policy.

We also observe that the dynamic programming framework allows nodes to find their optimal forwarding policy based on the statistics of their outgoing links and the "offered deadline-constrained utilities'' $U_{j}^{\star}(t+1)$ of one-hop away parents in a distributed fashion, see Algorithm~\ref{alg2:MultiPathRoutingDAG}. This step can be done a priori if the link statistics do not change, or at a slow time-scale using message-passing between nodes.
%
The complexity of the dynamic programming framework is of order ${\mathcal O}(D \vert \mathcal N\vert^2)$: the $\vert {\mathcal N}\vert $ nodes have to be scanned from the destination towards the source and from the deadline $t=D-1$ backwards in time to  $t=0$. In each such iteration, each node has to consider the option of forwarding to each of its maximum $\vert {\mathcal N}\vert$ outgoing neighbors.

\begin{algorithm}
\caption{Distributed implementation of the DP}
\begin{algorithmic}[t]
\label{alg2:MultiPathRoutingDAG}
    \STATE $U^\star_i(D)=0 \; \forall i \neq
Z
    ,$
    \quad $U^\star_
Z
    (t)=1 \; \forall t \in [0,D]$\\
    \FOR{$t=D-1$ to $0$}
        \FOR{$i=1$ to $
Z
        -1$}
            \STATE Compute maximum utility $U^\star_i(t)$
             by Eq.~\eqref{eq:MaximumUtility},
            \STATE Compute optimal action $j^\star_i(t)$
             by Eq.~\eqref{eq:OptimalPolicy},
            \STATE Transmit $U^\star_i(t)$ reliably to neighbor $j$ where $i \in \mathcal{N}_j$.
        \ENDFOR
    \ENDFOR
\end{algorithmic}
\end{algorithm}

\subsection{Linear-quadratic Gaussian control for fixed forwarding policy}\label{subsec:control_design}
Under a deadline-constrained forwarding policy $\pi$, the network delivers sensor packets with a fixed delay of $\tau=Dt_s$ seconds, and loses samples independently with probability $1-\rho^{\pi}$. In what follows, we will drop superscript $\pi$ for simplicity of notation. Moreover, we denote the process state by $x_{k}=x(kh)$, the control signal by $u_{k}=u(kh)$, the process noise by $v_{k}=v(kh)$, and measurement noise by $w_{k}=w(kh)$.

The evolution of the system between sampling instants can be described in terms of the extended state vector $\xi_{k}\triangleq\mbox{col}\{x_{k},u_{k-1}\}$ as
\begin{align*}
\xi_{k+1} = &\; \underbrace{\begin{bmatrix} e^{Ah} & \int_{h-\tau}^{h}e^{As}dsB \\ 0_{m\times n} & 0_{m\times m} \end{bmatrix}}_{\triangleq\; \Phi(h,\tau)} \xi_{k} + \underbrace{\begin{bmatrix} \int_{0}^{h-\tau}e^{As}dsB \\ I_{m\times m} \end{bmatrix}}_{\triangleq\; \Gamma(h,\tau)} u_{k} + \underbrace{\begin{bmatrix} I_{n\times n} \\ 0_{m\times n} \end{bmatrix}}_{\triangleq\; G} v_{k} \;,\\
y_k = &\;
\rho_{k}\underbrace{\begin{bmatrix} \tilde{C} & 0_{1\times m} \end{bmatrix}}_{\triangleq\; C}\xi_{k} + w_{k}\;,
\end{align*}
where $v_k$ and $w_k$ are zero mean discrete-time Gaussian white noise process with
\begin{equation}
\mathbf{E}\Bigg\{\begin{bmatrix} v_{k} \\ w_{k}\end{bmatrix} \begin{bmatrix} v_{k}^{\intercal} & w_{k}^{\intercal}\end{bmatrix} \Bigg\}
= \begin{bmatrix} R_{v} & 0 \\ 0 & R_{w} \end{bmatrix}\;,
\end{equation}
where $R_{v}\triangleq\int_{0}^{h}e^{As}R^{c}_{v}e^{A^{\intercal}s}ds$.

The continuous-time loss function~\eqref{eqn:continuous_loss} can be transformed into an equivalent discrete-time loss
\begin{align*}
\renewcommand{\arraystretch}{1.00}
\setlength{\arraycolsep}{2.00pt}
J_{d} = \mathbf{E}\Bigg\{\displaystyle{\sum_{k=0}^{N-1}}
\left[ \begin{array}{c}
\xi_{k} \\ u_{k}
\end{array} \right]^{\intercal}
\left[ \begin{array}{cc}
\Xi_{\xi\xi}(h,\tau) & \Xi_{\xi u}(h,\tau) \\
\Xi_{\xi u}^{\intercal}(h,\tau) & \Xi_{uu}(h,\tau)
\end{array} \right]
\left[ \begin{array}{c}
\xi_{k} \\ u_{k}
\end{array} \right] + \xi^{\intercal}_{N}\Xi_{0}\xi_{N}^{}\Bigg\} \;,
\end{align*}
where $N=\lceil\frac{\mathrm{T}}{h}\rceil$,
\begin{align*}
\Xi_{\xi\xi}(h,\tau) \triangleq &
\begin{bmatrix}
Q_{xx}^{\tau}+\Phi^{\intercal}(\tau)Q_{xx}^{h-\tau}\Phi(\tau) & Q_{xu}^{\tau}+\Phi^{\intercal}(\tau)Q_{xx}^{h-\tau}\Gamma(\tau) \\
Q_{xu}^{\tau\intercal}+\Gamma^{\intercal}(\tau)Q_{xx}^{h-\tau}\Phi(\tau)   & Q_{uu}^{\tau}+\Gamma^{\intercal}(\tau)Q_{xx}^{h-\tau}\Gamma(\tau)
\end{bmatrix}, \\
\Xi_{\xi u}(h,\tau) \triangleq &
\begin{bmatrix}
\Phi^{\intercal}(\tau)Q_{xu}^{h-\tau} \\
\Gamma^{\intercal}(\tau)Q_{xu}^{h-\tau}
\end{bmatrix}, \\
\Xi_{uu}(h,\tau) \triangleq &\; Q_{uu}^{h-\tau}
\end{align*}
\noindent with $\Phi(t)=e^{At}$, $\Gamma(t)=\int_{0}^{t}\Phi(s)Bds$, $Q_{xx}^{t}=\int_{0}^{t}\Phi^{T}(s)Q^{c}_{xx}\Phi(s)ds$,
$Q_{xu}^{t}=\int_{0}^{t}\Phi^{T}(s)\big(Q^{c}_{xx}\Gamma(s)+Q^{c}_{xu}\big)ds$ and $Q_{uu}^{t}=\int_{0}^{t}\big(\Gamma^{T}(s)Q^{c}_{xx}\Gamma(s)+2\Gamma^{T}(s)Q^{c}_{xu}+Q^{c}_{uu}\big)ds$.

The optimal control problem is then to compute the control sequence $\{u_k\}_{k\geq 0}$ that minimizes the discrete-time loss function. Note that $u_k$ is not computed until time $t=kh+\tau$, at which time $y_k$ is available to the controller unless it has been dropped by the network. Hence, the controller has access to the following information set when computing $u_k$:
\begin{equation*}
\mathcal{I}_{k} \triangleq
\big\{\mathcal{Y}_{k}\,,\,\mathcal{U}_{k-1}\,,\,\mathcal{R}_{k}\big\} \,.
\end{equation*}
Here, $\mathcal{Y}_{k}=\big(y_{k},\dots,y_{1}\big)$, and $\mathcal{U}_{k-1}=\big(u_{k-1},\dots,u_{1}\big)$, while $\mathcal{R}_{k}=\big(\rho_{k},\dots,\rho_{1}\big)$ is the
realizations of the Bernoulli random variable $\rho_k$ that models successful packet transmissions.
It is important to note that the discrete-time loss has cross-terms even if the continuous-time loss function does not. Schenato~\textit{et al.}~\cite{SSFPS:07} studied a similar problem without cross-terms in the loss function. In what follows, we extend the framework of~\cite{SSFPS:07} to include the cross-coupling terms in the loss function and derive the optimal controller and bound its performance.

\subsubsection{Estimator Design}
As in~\cite{SSFPS:07} the Kalman filter is the optimal estimator for our setting. The minimum mean square error (MMSE) estimate $\hat{\xi}_{k|k}$ of $\xi_{k}$ given by $\hat{\xi}_{k|k}=\mathbf{E}\{\xi_{k}|\;\mathcal{I}_{k}\}$ can be computed recursively starting from the initial conditions $\hat{\xi}_{0|-1}=\mbox{col}\{0_{n\times 1},~0_{m\times 1}\}$ and $P_{0|-1}=P_{0}$. The innovation step is
\begin{align}
\hat{\xi}_{k+1|k}&\triangleq\mathbf{E}\{\xi_{k+1}|\;\mathcal{I}_{k}\}
=\Phi\hat{\xi}_{k|k}+\Gamma u_{k}\label{eqn:kalman_first}\\
e_{k+1|k}&\triangleq \xi_{k+1}-\hat{\xi}_{k+1\vert k}=\Phi e_{k|k} + Gv_{k}
\\
P_{k+1|k}&\triangleq\mathbf{E}\{e_{k+1|k}e^{\intercal}_{k+1|k}|\;
\mathcal{I}_{k}\}=\Phi P_{k|k}\Phi^{\intercal}+\tilde{R}_{v}\label{eq:P_first}
\end{align}
where $\tilde{R}_{v}\triangleq GR_{v}G^{\intercal}$ and $v_{k}$ is independent from $\mathcal{I}_{k}$, while the correction step is
\begin{align}
\hat{\xi}_{k+1\vert k+1}&=\hat{\xi}_{k+1\vert k}+\rho_{k+1}K_{k+1}(y_{k+1}-C\hat{\xi}_{k+1\vert k})\\
e_{k+1\vert k+1}&=\xi_{k+1}-\hat{\xi}_{k+1\vert k+1}\\
P_{k+1\vert k+1}&=P_{k+1\vert k}-\rho_{k+1}K_{k+1}CP_{k+1\vert k}\label{eq:P_second}\\
K_{k+1}&\triangleq
P_{k+1\vert k}C^{\intercal}(CP_{k+1\vert k}C^{\intercal}+R_{w})^{-1}\,.\label{eq:kalman_gain}
\end{align}
The following result, similar to~\cite{SSFPS:07}, characterizes the estimation error covariance matrix.
\begin{proposition} \label{prop:kalman} The MMSE estimate $\hat{\xi}_{k\vert k}$ of $\xi_k$ is given by the time-varying Kalman filter (\ref{eqn:kalman_first}) -- (\ref{eq:kalman_gain}). The expected value of the covariance matrix can be bounded as
\begin{align*}
\underline{P}_{k\vert k}&\leq \mathbf{E}_{\rho}\{ P_{k\vert k}\}\leq \overline{P}_{k \vert k}
\end{align*}
where the bounds can be computed iteratively as
\begin{align*}
\renewcommand{\arraystretch}{1.4}
\setlength{\arraycolsep}{1.6pt}
\begin{array}{rcl}
\overline{P}_{k+1|k} & = &\Phi\overline{P}_{k|k-1}\Phi^{\intercal} +\tilde{R}_{v} - \rho\Phi\overline{P}_{k|k-1}C^{\intercal}(C\overline{P}_{k|k-1}C^{\intercal}+R_{w})^{-1}C\overline{P}_{k|k-1}\Phi^{\intercal} \\
\overline{P}_{k|k} & = & \overline{P}_{k|k-1}-\rho\overline{P}_{k|k-1}C^{\intercal}(C\overline{P}_{k|k-1}C^{\intercal}+R_{w})^{-1}C\overline{P}_{k|k-1} \\
\underline{P}_{k+1|k} & = & (1-\rho)\Phi\underline{P}_{k|k-1}\Phi^{\intercal}+\tilde{R}_{v} \\
\underline{P}_{k|k} & = & (1-\rho)\underline{P}_{k|k-1}
\end{array}
\end{align*}
starting from  the initial conditions $\underline{P}_{0|-1}=\overline{P}_{0|-1}=P_{0}$. When $k\rightarrow \infty$, the iterations converge to the unique stationary solutions $\underline{P}_{\infty}$ and $\overline{P}_{\infty}$ of the modified algebraic Riccati equations
\begin{align}
\overline{P}_{\infty} & = \Phi\overline{P}_{\infty}\Phi^{\intercal}+\tilde{R}_{v}-\rho\Phi\overline{P}_{\infty}C^{\intercal}(C\overline{P}_{\infty}C^{\intercal}+R_{w})^{-1}C\overline{P}_{\infty}\Phi^{\intercal}\;, \label{eqn:mare1} \\
\underline{P}_{\infty} & = (1-\rho)\Phi\underline{P}_{\infty}\Phi^{\intercal}+\tilde{R}_{v}\;. \label{eqn:mare2}
\end{align}
\end{proposition}

\subsubsection{Controller Design}
Next, we develop the optimal state feedback control law.

\begin{proposition}\label{prop:control}
Consider the finite horizon LQG control problem. The optimal control law
\begin{equation}
u_{k}=-\underbrace{(\Gamma^{\intercal}S_{k+1}\Gamma+\Xi_{uu})^{-1}(\Gamma^{\intercal}S_{k+1}\Phi+\Xi^{\intercal}_{\xi u})}_{\triangleq L_{k}}\hat{\xi}_{k|k}\label{eq:static_controller}
\end{equation}
is a linear function of the estimated state. The matrix $S_{k}$ evolves according to the recursion
\begin{align}
S_{k}=\;\Phi^{\intercal}S_{k+1}\Phi +\Xi_{\xi\xi} - (\Phi^{\intercal}S_{k+1}\Gamma+\Xi_{\xi u})(\Gamma^{\intercal}S_{k+1}\Gamma +\Xi_{uu})^{-1}(\Gamma^{\intercal}S_{k+1}\Phi+\Xi^{\intercal}_{\xi u})\;, \label{eqn:control_riccati}
\end{align}
where $\hat{\xi}_{k|k}$ is the MMSE estimate of the state $\xi_{k}$ based on the information set $\mathcal{I}_k$ computed with the Kalman filter (\ref{eqn:kalman_first})--(\ref{eq:kalman_gain}). As $k\rightarrow\infty$, the recursion converges to a unique solution
\begin{equation*}
S_{\infty} =\;\Phi^{\intercal}S_{\infty}\Phi+\Xi_{\xi\xi} - (\Phi^{\intercal}S_{\infty}\Gamma+\Xi_{\xi u}) (\Gamma^{\intercal}S_{\infty}\Gamma+\Xi_{uu})^{-1}(\Gamma^{\intercal}S_{\infty}\Phi+\Xi^{\intercal}_{\xi u}) \;,
\end{equation*}
for which the associated stationary controller gain is
\begin{align*}
L_{\infty}\triangleq \lim_{k\rightarrow\infty}L_{k}=-(\Gamma^{\intercal}S_{\infty}\Gamma+\Xi_{uu})^{-1}(\Gamma^{\intercal}S_{\infty}\Phi+\Xi_{\xi u}^{\intercal})\;.
\end{align*}
\end{proposition}
{\noindent \textit{Proof:}}
To derive the optimal feedback control law and the corresponding value for the objective function, we apply dynamic programming. We let the optimal value function $V_{k}(\xi_{k})$ be
\begin{align}
V_{k}(\xi_{k})\triangleq &\;\min_{u_{k}}\mathbf{E}\{\xi_{k}^{\intercal}\Xi_{\xi\xi}\xi_{k}+2\xi_{k}^{\intercal}\Xi_{\xi u}u_{k}+u_{k}^{\intercal}\Xi_{uu}u_{k}+V_{k+1}\vert~\mathcal{I}_{k}\},\label{eq:Bellman_func}\\
V_{N}(\xi_{N})\triangleq &\;\mathbf{E}\{\xi_{N}^{\intercal}\Xi_{0}\xi_{N}|~\mathcal{I}_{N}\},\label{eq:Bellman_func_initial}
\end{align}
where $k=\{N-1,\ldots,1\}$. We show that $J_{N}^{\star}=V_{0}(\xi_{0})$. The solution of the Bellman equation~\eqref{eq:Bellman_func} with the initial condition~\eqref{eq:Bellman_func_initial} is given by
\begin{align}
V_{k}(\xi_{k})=\mathbf{E}\{\xi_{k}^{\intercal}S_{k}\xi_{k}|\;\mathcal{I}_{k}\}+c_{k}\,,\label{eq:Bellman_func_solution}
\end{align}
where the nonnegative matrix $S_{k}$ and the scalar $c_{k}$ are independent of the information set $\mathcal{I}_{k}$. In contrast to~\cite{SSFPS:07}, we allow for cross-coupling terms in the Bellman equation~\eqref{eq:Bellman_func}, which is critical for comparing~the control costs under different sampling intervals. Apart from this, the proof is similar to the one in~\cite{AW:96,SSFPS:07}.
\hfill $\blacksquare$

\subsubsection{Optimal control cost}
The optimal loss function of the finite horizon LQG problem is
\begin{multline}
J^{\star}_{N}(\rho)=\xi^{\intercal}_{0}S_{0}\xi_{0}+\mbox{Tr}\big(S_{0}P_{0}\big)+\sum_{k=0}^{N-1}\mbox{Tr}\big(S_{k+1}\tilde{R}_{v}\big) \\
  \hspace{10mm}   +\sum_{k=0}^{N-1}\mbox{Tr}\big((\Phi^{\intercal}S_{k+1}\Phi+\Xi_{\xi\xi}-S_{k})\mathbf{E}_{\rho}\{P_{k|k}\}\big), \label{eq:network_loss_func}
\end{multline} %
where expectation is taken over a Bernoulli sequence $\{\rho_k\}_{k\geq 0}$ with $\mathbf{E}\{\rho_k\}=\rho$. Since no efficient way of computing the expectation is known, one can use the upper and lower bounds on $\mathbf{E}_{\rho}\{P_{k\vert k}\}$ given in Proposition~\ref{prop:kalman} to compute the upper and lower bounds on the finite-horizon control cost $J_N^{\min}(\rho)\leq J_N^{\star}(\rho)\leq J_N^{\max}(\rho)$. For the infinite horizon case, the bounds become
\begin{align}
J_{\infty}^{\min} \triangleq &\; \lim_{N\rightarrow\infty}\frac{1}{N}J^{\min}_{N} \nonumber\\
 =  &\; \mbox{Tr}(S_{\infty}\tilde{R}_{v})+(1-\rho)\mbox{Tr}\big((\Phi^{\intercal}S_{\infty}\Phi+\Xi_{\xi\xi}-S_{\infty})\underline{P}_{\infty}\big),
\label{eqn:min_bound_loss_function}\\
J_{\infty}^{\max} \triangleq &\; \lim_{N\rightarrow\infty}\frac{1}{N}J^{\max}_{N} \nonumber\\
 = &\; \mbox{Tr}(S_{\infty}\tilde{R}_{v}) + \mbox{Tr}\big((\Phi^{\intercal}S_{\infty}\Phi+\Xi_{\xi\xi}-S_{\infty})(\overline{P}_{\infty}-\rho\overline{P}_{\infty}C^{\intercal} \nonumber\\
&\; \times(C\overline{P}_{\infty}C^{\intercal}+R_{w})^{-1}C\overline{P}_{\infty})\big), \label{eqn:max_bound_loss_function}
\end{align}
where the matrices $\underline{P}_{\infty}$, $\overline{P}_{\infty}$ and $S_{\infty}$, are given in Proposition~\ref{prop:kalman} and Proposition~\ref{prop:control}, respectively.

In brief, the optimal estimator is a time-varying Kalman filter given by~\eqref{eq:P_first},~\eqref{eq:P_second} and~\eqref{eq:kalman_gain}, while the control law is a static linear feedback~\eqref{eq:static_controller}. The combined performance, in the sense of the continuous-time loss function~\eqref{eqn:continuous_loss}, can be bounded as
in~(\ref{eqn:min_bound_loss_function}) and~(\ref{eqn:max_bound_loss_function}). It is this controller and these performance bounds that we use in our co-design.

\subsection{Optimality of the co-design framework} \label{sec:optimality}

Next, we provide a formal proof that the achievable loss for a fixed sampling interval and time lag is increasing in the loss probability, which allows us to conclude that the combination of linear-quadratic control under loss and deadline-constrained maximum reliability routing provides the optimal co-design.

\begin{lemma}\label{lem:mono2}
Consider operators
$f(X)=\Phi X \Phi^{\intercal}+\tilde{R}_{v}$, $h_{\rho}(X)=X-\rho XC^{\intercal}(CXC^{\intercal}+R_w)^{-1}CX$ and $g_{\rho}(X)=h_{\rho}(f(X))$. If $X^{\prime}\leq X$  and $\rho^{\prime}\geq \rho$, then $g_{\rho^{\prime}}(X^{\prime}) \leq g_{\rho}(X)$.
\end{lemma}
{\noindent \textit{Proof:}}
By Lemma~\ref{lem:mono}~(a,c), $g_{\rho^{\prime}}(X^{\prime})\leq g_{\rho}(X^{\prime})\leq g_{\rho}(X)$.
\hfill $\blacksquare$

\begin{theorem} \label{thm:monotonicity}
For given sampling interval $h$ and time lag $\tau$, the optimal control loss $J_{N}^{\star}(\rho)$ is monotone decreasing in end-to-end reliability $\rho$.
\end{theorem}
{\noindent \textit{Proof:}}
Our proof relies on a coupling argument~\cite{Lin:92} on the underlying end-to-end loss processes. Specifically, for any two Bernoulli processes $\{ \rho_k \}_{k\geq 0}$ and $\{ \rho^{\prime}_k\}_{k\geq 0}$ with $\mathbf{E}\{\rho_k\} = \rho$ and $\mathbf{E}\{\rho^{\prime}_k\} = \rho^{\prime}$ such that $\rho\leq \rho^{\prime}$ we establish that  $J_N^{\star}(\rho^{\prime})\leq J_N^{\star}(\rho)$.

The key idea of coupling is to define the two processes on a common probability space on which the analysis is carried out. To this end, let $\{\omega_k\}_{k\geq 0}$ be a sequence of independent random variables with $\omega_k\thicksim \mathrm{Uni}(0,1)$ and define $\rho_k = \mathds{1}_{\omega_k \leq \rho}$ and $\rho^{\prime}_k=\mathds{1}_{\omega_k\leq \rho^{\prime}}$. Then, $\{\rho_k\}_{k\geq 0}$ and $\{\rho_k^{\prime}\}_{k\geq 0}$ are Bernoulli trial processes with probabilities $\rho$ and $\rho^{\prime}$, respectively. A crucial property of our construction is that for every realization $\{\omega_k\}_{k\geq 0}$, the associated sequences $\{\rho_k(\omega_k)\}_{k\geq 0}$ and $\{\rho_k^{\prime}(\omega_k)\}_{k\geq 0}$ satisfy $\rho_k^{\prime}(\omega_k)\geq \rho_k(\omega_k)$ for all $k\in\mathbb{N}$. In particular, whenever the more reliable sequence has a loss, the less reliable sequence also has a loss.

Our next step is to show that $P_{N\vert N}^{\prime}\leq P_{N\vert N}$ where $P_{N\vert N}^{\prime}$ is the estimation error covariance matrix under packet delivery sequence $\{ \rho_k\}_{k\geq 0}$ and initial value $P_{0\vert 0}^{\prime}=P_0$ and $P_{N\vert N}$ is the estimation error covariance matrix under packet delivery sequence $\{\rho_k^{\prime}\}_{k\geq 0}$ and the same initial value $P_{0\vert 0}=P_0$. We will establish this claim using induction. Clearly, $P^{\prime}_{k\vert k}\leq P_{k\vert k}$ holds for $k=0$. Assume that it holds for an arbitrary $k$. Then, since $P_{k+1\vert k+1}=g_{\rho}(P_{k\vert k})$ and $P_{k+1\vert k+1}^{\prime}=g_{\rho^{\prime}}(P^{\prime}_{k\vert k})$, Lemma~\ref{lem:mono2} implies that $P_{k+1\vert k+1}^{\prime}\leq P_{k+1\vert k+1}$. Hence, by induction, $P_{N\vert N}^{\prime}\leq P_{N\vert N}$.

Finally, combining this observation with Lemma~\ref{lem:trace} (in Appendix):
\begin{align*}
J_N^{\star}(\rho) &= c + \underset{\rho_1, \dots, \rho_N}{\mathbf{E}}\, \sum_{k=0}^{N-1} \mbox{Tr}\, \Delta_k P_{k\vert k}(\rho) \\
& = c+ \underset{\omega_1, \dots, \omega_N}{\mathbf E}\, \sum_{k=0}^{N-1} \mbox{Tr}\, \Delta_k P_{k\vert k}(\rho(\omega_{k})) \\
&\leq c+\underset{\omega_1, \dots, \omega_N}{\mathbf E}\,\sum_{k=0}^{N-1} \mbox{Tr}\, \Delta_k P_{k\vert k}(\rho^{\prime}(\omega_{k})) \\
& = c+\underset{\rho_1^{\prime}, \dots, \rho_N^{\prime}}{\mathbf E} \sum_{k=0}^{N-1} \mbox{Tr}\, \Delta_k P_{k\vert k}(\rho^{\prime})  =J_N^{\star}(\rho^{\prime})\;.
\end{align*}
This concludes the proof.
\hfill $\blacksquare$

Theorem~\ref{thm:monotonicity} allows us to establish optimality of our co-design framework.
\begin{theorem}\label{thm:codesign_optimality}
For given sampling interval and time lag, the optimal closed loop performance, in the sense of the linear-quadratic control loss (\ref{eqn:continuous_loss}), is obtained by scheduling the network to maximize deadline-constrained reliability and computing the control action by the time-varying LQG-controller (\ref{eqn:kalman_first})-(\ref{eq:kalman_gain}), (\ref{eq:static_controller}), (\ref{eqn:control_riccati}).
\end{theorem}

Theorem~\ref{thm:codesign_optimality} establishes that the optimal co-design can be obtained as follows:
we sweep over both the sampling interval $h$ and the time lag $\tau$, schedule the network to maximize on-time delivery probability,  and estimate the associated closed-loop performance of the time-varying LQG controller for the corresponding end-to-end reliability $\rho$. The optimal co-design is obtained for the sampling interval and the time lag that attain the minimal closed-loop loss.

It is natural to ask if one could find the optimal $h$ and $\tau$ in a more efficient way than by exhaustive search. We will show in the numerical examples that the optimal loss as a function of $h$ might exhibit multiple local minima, which indicates that it would be challenging to find a universal and efficient way of picking the optimal sampling interval. Another issue is that there is no efficient way to compute the  expectation of the covariance matrix with respect to the loss process apart from \emph{e.g.}, Monte Carlo simulation~\cite{SSF+:04}. To overcome this problem, one could replace the true performance expression by the upper bound $J_N^{\max}(\rho)$ and pick the sampling interval that minimizes this upper bound. The next result establishes that $J_N^{\max}$ has the same monotonicity properties as $J_N^{\star}$, and we will show in the numerical examples that $J_N^{\max}$ provides a good surrogate to the true loss function when it comes to picking the optimal sampling interval.
\begin{lemma}
The upper bound on the control loss $J_{N}^{\max}(\rho)$ is monotone decreasing in end-to-end reliability $\rho$.
\end{lemma}
{\noindent \textit{Proof:}} Note $J_{\infty}^{\max}(\rho) = c^{\prime}+\mbox{Tr}\, \Delta_{\infty}h_{\rho}(\overline{P}_{\infty}(\rho))$ where $\overline{P}_{\infty}(\rho)$ is the stationary solution of the modified algebraic Riccati equation (\ref{eqn:mare1}) for reliability parameter $\rho$.  Now, consider two reliability parameters $\rho$ and $\rho^{\prime}$. By similar arguments in Theorem~\ref{thm:monotonicity}, $\rho^{\prime}\geq \rho$ implies that $\overline{P}_{\infty}(\rho^{\prime})\leq \overline{P}_{\infty}(\rho)$ and by Lemma~\ref{lem:mono}(b), $h_{\rho}(\overline{P}_{\infty}(\rho))\leq h_{\rho^{\prime}}(\overline{P}_{\infty}(\rho^{\prime}))$. By Lemma~\ref{lem:trace}, $J_{\infty}^{\max}(\rho^{\prime})\leq J_{\infty}^{\max}(\rho)$. Our claim is proved.
\hfill $\blacksquare$

\section{Numerical Examples}\label{sec:Numerical_Examples}
\begin{figure}[!htbp]
\centerline{\includegraphics[width=0.45\columnwidth]
  {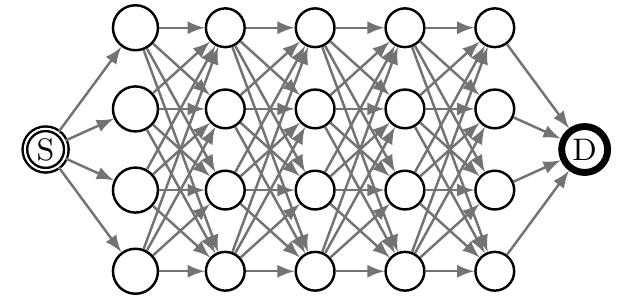}}
  \caption{\normalsize Network topology with the source 6-hop from the destination.}
  \label{fig:DODAG}
\end{figure}

We now demonstrate our co-design procedure on numerical examples. We consider the following second-order linear system
\begin{align}
\renewcommand{\arraystretch}{1.00}
\setlength{\arraycolsep}{2.25pt}
\begin{array}{rcl}
dx & = & \left[\begin{array}{cc} 0 & 1 \\ -\omega_{0}^{2} & -2\alpha\zeta\omega_{0} \end{array}\right]xdt + \left[\begin{array}{c} 0 \\ \omega_{0}^{2} \end{array}\right]udt +dv_{c} \;, \\
y(kh) & = & \left[%
\begin{array}{cc}
  1 & 0 \\
\end{array}%
\right]x(kh)+w(kh),
\end{array}
\label{eq:num_example}
\end{align}
where $v_{c}$ has incremental covariance $R^{c}_{v}= \mbox{diag}\{0.5,0.5\}$, and $w(kh)$ has covariance $R_{w}=10^{-4}$. Our co-design should minimize~\eqref{eqn:continuous_loss} for $Q_{xx}^{c}=\mbox{diag}\{2,1\}$, $Q_{xu}^{c}=0_{2\times 1}$, $Q_{uu}^{c}=1$, and $Q_{0}^{c}=0_{2\times 2}$. The multi-hop wireless network between sensor and controller is shown in Fig~\ref{fig:DODAG}. A link between two nodes exists if and only if the transmission loss probability between these two nodes is strictly smaller than $1$ (\emph{i.e.}, if communication between the two nodes is at all possible). The length of one time slot is $t_s=10$ms. Furthermore, in all numerical examples, we assume that the time lag $\tau$ is always set to equal the sampling period $h$.

\subsection{No energy constraint}
We here consider the case without energy constraints in which the optimal forwarding policy and the maximum deadline-constrained reliability are derived by the dynamic programming framework with $\delta=0$ in \S~\ref{sec:DPSolution}.

\begin{figure}[!htbp]
\centerline{\includegraphics[width=0.6\columnwidth]
  {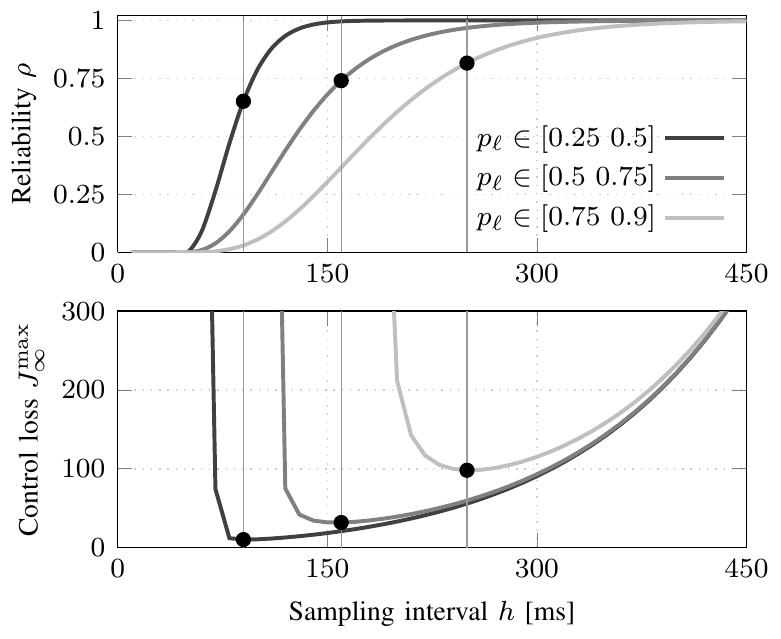}}
  \caption{\normalsize Comparison of the upper bounds $J^{\max}_{\infty}$ for three different network scenarios. Note that $p_{\ell}$ is the link loss probability.}
  \label{fig:case_study_infinite}
\end{figure}

We firstly consider the system~\eqref{eq:num_example} to be unstable with the parameters $\alpha=-1$, $\zeta=1$, and $\omega_{0}=1$. Periodic samples of the output $y(kh)$ are transmitted over the network in Fig.~\ref{fig:DODAG}. For every sampling period $h$, we compute the optimal forwarding policy and the associated deadline-constrained reliability. The achievable latency-reliability pairs for three network scenarios in which the network becomes increasingly unreliable (by increasing the loss probabilities on links) are shown in Fig.~\ref{fig:case_study_infinite} (top). We then discretize the control loss function~\eqref{eqn:continuous_loss}, solve the corresponding LQG optimal control problem, and evaluate its performance as described in~\S~\ref{subsec:control_design}.
Fig.~\ref{fig:case_study_infinite} (bottom) shows the optimal closed-loop control losses for varying sampling intervals under the three network scenarios. We note that there is no single optimal sampling interval or target end-to-end reliability. The optimal sampling interval ranges from $90\mathrm{ms}$, for the most reliable network scenario, to $250\mathrm{ms}$, for the least reliable case. This corresponds to a required end-to-end reliability of $65\%$ and $82\%$, respectively. As the network becomes less reliable,  more retransmissions (a longer $h^\star$) are required to guarantee a sufficiently high reliability which causes the associated control loss to increase. Specifically, the optimal control loss (marked with a dot in Fig.~\ref{fig:case_study_infinite}) increases by approximately a factor of ten from the most to the least reliable network scenarios.
%
\begin{table*}[!htbp]
\renewcommand{\arraystretch}{1.2}
\centering
\begin{tabular}
{
 >{\centering \arraybackslash \columncolor{black!05!white}} c <{}
 >{\centering \arraybackslash } c <{}
 >{\centering \arraybackslash \columncolor{black!05!white}} c <{}
 >{\centering \arraybackslash } c <{}
}
 \toprule
 \bfseries $\zeta$ &
 \bfseries $\omega_{0}$ &
 $h$ (rule-of-thumb) &
 $h^{\star}$ (Optimal sampling interval) \\
 \midrule
 0.1 &
 \multicolumn{1}{c}{\multirow{3}{*}{$\frac{\pi}{\sqrt{1-\zeta^{2}}}$}} &
 40 \; -- \;90~ms &
 130~ms \\
 0.2 &
 \multicolumn{1}{c}{} &
 40 \; -- \;100~ms &
 120~ms \\
 0.7 &
 \multicolumn{1}{c}{} &
 50 \;-- \;120~ms &
 60~ms \\
 \bottomrule
\end{tabular}
\caption{Selection of sampling interval $h$ in the wireless control system that includes communication links with link loss probability $p_{\ell}\in[0.75,0.9].$} \label{table_example}
\end{table*}

\begin{figure}[!htbp]
\centerline{\includegraphics[width=0.6\columnwidth]
  {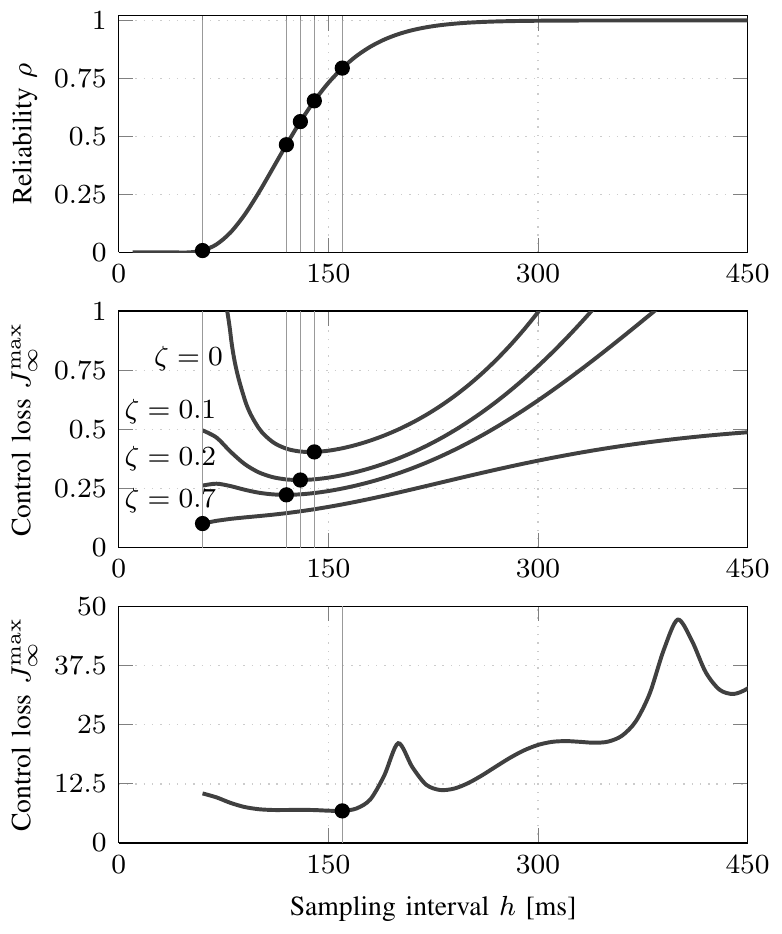}}
  \caption{\normalsize Comparison of the upper bounds $J_{\infty}^{\max}$ for the system~\eqref{eq:num_example} which is parameterized by several damping ratios $\zeta$ and natural frequencies $\omega_{0}$.}
 \label{fig:experimental_delay}
\end{figure}

We now consider the case when the system~\eqref{eq:num_example} is stable, and $\zeta\geq 0$, $\alpha=1$, and $\omega_{0}=\pi/\sqrt{1-\zeta^{2}}$. This parameterization is chosen roughly to yield the same optimal sampling interval independent of $\zeta$. Fig.~\ref{fig:experimental_delay} (center) shows how the trade-off curves change when we vary the damping of the system poles. Similar to the unstable open-loop systems, distinct co-design optima also exist for the open-loop stable systems, but they become less relevant when the system poles are close to being critically damped. Furthermore, we compare the optimal sampling interval obtained by our co-design framework with the rule-of-thumb from~\cite[pp.~129-130]{AW:96}. As can be seen in Table~\ref{table_example}, the optimal sampling interval is, more often than not, outside the range of sampling intervals proposed by the rule-of-thumb.

It may appear discomforting that our co-design framework requires sweeping over the sampling interval to find the jointly optimal design. However, as shown in Fig.~\ref{fig:experimental_delay} (bottom) with parameters $\alpha=1$, $\zeta=0.02$ and $\omega_{0}=5\pi$, the optimal closed-loop control cost as a function of the sampling interval might exhibit multiple local minima, which indicates that it will be hard to circumvent this search in general.

In this section, we set the maximum transmission delay equal to the sampling interval to reduce the search space to a single parameter. One could imagine that a better performance could be achieved by optimizing over both the sampling interval and the transmission delay, which is indeed possible in our framework. While we cannot rule out this possibility in general, our experience is that the additional complexity of a two-dimensional search does not drastically improve the performance.

\begin{figure}[!htbp]
\centerline{\includegraphics[width=0.6\columnwidth]
  {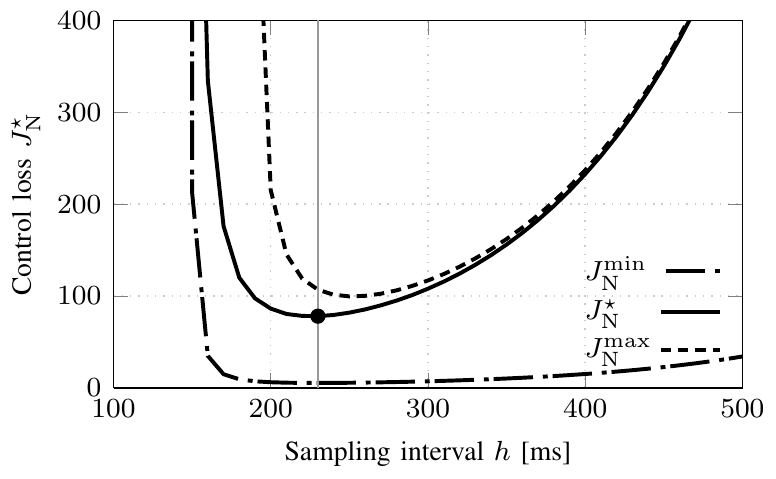}}
  \caption{\normalsize The experimental $J^{\star}_{N}$, upper bound $J^{\max}_{N}$ and lower bound $J^{\min}_{N}$ for the minimum control cost with respect to sampling period. The curves are obtained by averaging 10,000 Monte Carlo simulations for the horizon length $T=500\mathrm{s}$, with the arrival sequence $\{\rho_{k}\}_{k\geq 0}^{}$ generated randomly.}
 \label{fig:case_study_finite}
\end{figure}

Lastly, since the true closed-loop performance can only be evaluated using Monte Carlo simulations, it is useful to be able to rely on more easily computable upper and lower performance bounds, introduced in this paper, in the search for the optimal sampling interval. Fig.~\ref{fig:case_study_finite} compares the upper, lower, and true closed-loop performance for the co-design under the least reliable network scenario shown in Fig.~\ref{fig:experimental_delay} (top), and the system~\eqref{eq:num_example} with parameters $\alpha=-1$, $\zeta=1$, and $\omega_{0}=1$. The upper bound $J_N^{\max}$ becomes quite accurate for sampling intervals when $h\geq250\mathrm{ms}$, and the sampling interval that minimizes the upper bound on the performance is rather close to the optimal sampling interval for the true cost. Hence, we believe that the upper bound is a good surrogate for the true performance if we need to carry out the co-design with limited computational resources.

\subsection{With energy constraints}
\begin{figure}[t]
\centerline{\includegraphics[width=0.65\columnwidth]
  {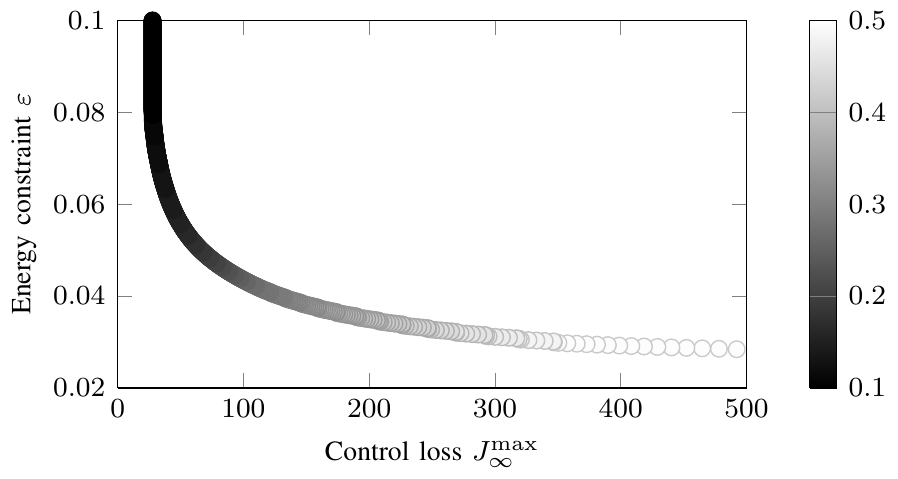}}
  \caption{\normalsize The optimal control loss for different energy cost constraints and the corresponding optimal sampling periods in seconds (shown in gray scale).}
 \label{fig:ControlVsEnergy2}
\end{figure}

We here include the energy constraint in the multi-hop wireless network, and consider the unstable system~\eqref{eq:num_example} with the parameters $\alpha=-1$, $\zeta=1$, and $\omega_{0}=1$. Fig.~\ref{fig:ControlVsEnergy2} illustrates the optimal closed-loop control losses for a range of energy cost constraints. Moreover, in the same plot, the optimal sampling intervals are shown in gray scale with shorter periods for smaller control losses. The optimal control performance naturally decreases when the energy cost constraint becomes increasingly stringent. An intriguing observation is that it is exceedingly costly to achieve the minimum control loss.  As revealed by Fig.~\ref{fig:ControlVsEnergy2}, significant energy savings can be ensured by accepting a relatively small deterioration in the control performance.

\begin{figure}[t]
\centerline{\includegraphics[width=0.55\columnwidth]
  {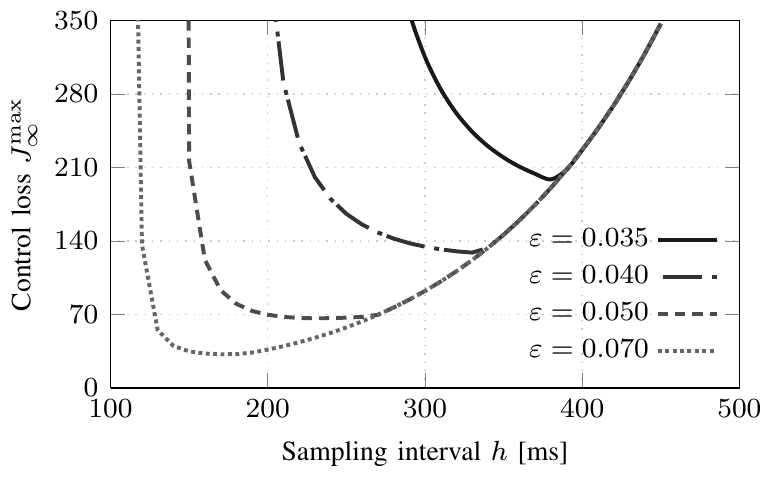}}
  \caption{\normalsize Comparison of the control loss upper bounds $J_{\infty}^{\max}$ with different sampling period for a set of energy constraints.}
 \label{fig:ControlVsEnergy1}
\end{figure}

Fig.~\ref{fig:ControlVsEnergy1} shows the control loss versus sampling periods for a set of energy constraints. Initially, with longer sampling periods, although fewer packets are injected into the network, these packets have larger deadline and more energy can be allocated to each packet with a higher per-packet delivery reliability. Hence, the control loss decreases. However, with even longer sampling periods, the energy constraint is no longer the bottleneck of the control performance, and the maximum reliability can be obtained without violating energy constraints. Thus, the control performance deteriorates since fewer packets are injected. This also explains why the curves overlap for large sampling periods. Finally, the phenomenon is more obvious with more stringent energy constraint where the system is more sensitive to different sampling intervals.

\section{Conclusions, Discussions and Future Works}\label{sec:conclusions}

We considered the joint design of forwarding policies and controllers for networked control loops that use multi-hop wireless communication for transmitting sensor data from process to controller. By parameterizing the design problem in terms of the sampling rate of the control loop, the co-design problem separates into two well-defined networking and control design tasks, both which admit optimal solutions: the network should be operated to maximize the deadline-constrained reliability subject to a total energy budget, and the control design should optimize closed-loop performance under packet loss. We develop optimal solutions to these problems, and demonstrate how the jointly optimal design can be found by a one-dimensional search over the sampling interval. To the best of our knowledge, this is the first co-design procedure that covers such a breadth of design parameters, gives a clean and modular design, yet guarantees to find the joint design with optimal closed-loop performance.

Admittedly, the framework has its limitations and several extensions are worthwhile to consider. One important extension is the multi-loop problem where multiple sensors take measurements and forward across the network. The corresponding networking task, however, becomes a real-time deadline-constrained multi-flow scheduling problem, which was recently proved to be NP-hard~\cite{Saifullah10}, so optimal designs will probably be hard to find. One interesting direction of research would be to develop sub-optimal policies for multiple real-time streams using, for example, approximate dynamic programming techniques. Another extension would be to consider more realistic physical models, such as nonlinear system dynamics and network models that allow for correlated link losses. When link losses are correlated, it is not the networking sub-problem that is major obstacle (see our paper~\cite{Zou2012} which develops real-time forwarding policies without energy constraints). It is the controller sub-problem that becomes hard, mainly since the end-to-end losses become correlated and the state of the Markov chain that describes end-to-end loss process is not available at the controller node. A third class of extensions would be to consider other control architectures. One interesting architecture is to stick to time-triggered sensing but allow for event-triggered control and actuator updates in the spirit of Nilsson~\cite{Nil:98}. The challenge of such architecture is that the system-level performance depends on the complete latency distribution of packets, which creates significant couplings between the controller and network designs. Another interesting controller architecture is to use event-triggered sampling, but the theory for event-triggered control is still in its infancy and many intermediate results need to fall in place before one could develop a corresponding co-design framework. Finally, one could also consider developing forwarding policies that require less network state information, use better models for energy consumption, and directly address the expected network lifetime. We hope to return to some of these extensions in our future work.

\paragraph*{Acknowledgements} The authors are grateful to Prof. Alexandre Proutiere for introducing them to the coupling method, and to Dr. Kin Sou and Dr. Daniel Lehmann for valuable comments on an earlier draft of this manuscript.

\section{Appendix}\label{sec:appex}

\label{sec:proofParetoFrontier}

\vspace{3mm}
{\noindent \textbf{Proof of Lemma~\ref{lem:ParetoFrontier}:}}
The number of actions in this MDP is limited by the number of neighboring nodes. We have a finite number of policies that lead to a finite number of energy costs. Thus, $\mathcal{C}$ is a finite~set.

For any given $\delta_1, \delta_2 \in \Delta_C$, we have $C^\star(\delta_1)=C^\star(\delta_2)$ with optimal polices $\pi^\star(\delta_1)$ and $\pi^\star(\delta_2)$, and reliability $\rho^\star(\delta_1)$ and $\rho^\star(\delta_2)$.
Without loss of generality, we let $\delta_1 < \delta_2$.
We can easily show that $\rho^\star(\delta)$ is a non-increasing function over $\delta \geq 0$, and hence $\rho^\star(\delta_1) \geq \rho^\star(\delta_2)$.
If $\rho^\star(\delta_1) = \rho^\star(\delta_2)$, then the lemma is proved. 
Now, we consider $\rho^\star(\delta_1) > \rho^\star(\delta_2)$.
Note the optimal utility with $\delta_2$ is $\rho^\star(\delta_2) - \delta_2 \cdot C^\star(\delta_2)$. However if we apply policy $\pi^\star(\delta_1)$, the utility with $\delta_2$ is $\rho^\star(\delta_1) - \delta_2 \cdot C^\star(\delta_1) > \rho^\star(\delta_2) - \delta_2 \cdot C^\star(\delta_2)$, which contradicts the optimality of $\pi^{\star}(\delta_2)$.
Hence, for any $C\in {\mathcal C}$, $\rho^\star(\delta)$ is unique for all $\delta \in \Delta_C$.
\hfill $\blacksquare$

\vspace{3mm}
{\noindent \textbf{Proof of Theorem~\ref{thm:ParetoFrontier}:}}
Let $g(\delta) \triangleq \max_\pi \{ \rho^\pi - \delta \cdot C^\pi \}
+ \delta \cdot C_\text{req}$ and $h(\delta)
\triangleq \max_\pi\{ \rho^\pi - \delta \cdot C^\pi \}$.
The proposed dynamic programming framework computes $h(\delta)$ for a given value of $\delta$ and returns a history-independent and deterministic policy. According to Markov decision process theory~\cite{Put:94}, $h(\delta)$ can be formulated as a linear program whose objective function coefficients depend on $\delta$, and it can be shown that $h(\delta)$ is a continuous function over $\delta$. Hence, $g(\delta)$ is also a continuous function over $\delta$.

It can be easily shown that $C^\star(\delta)$ is a non-increasing function over $\delta \geq 0$.
Hence, we have that $\Delta_{C^{(m)}}$ is an interval or a single point for a given $C^{(m)} \in \mathcal{C}$ by the fact $\mathcal{C}$ is a finite set from Lemma~\ref{lem:ParetoFrontier}.
Moreover, $h(\delta)=\rho^\star(\delta)-\delta C^\star(\delta)$ is a continuous function over $\delta$ and $\rho^\star(\delta)$ is unique for $\delta \in \Delta_{C^{(m)}}$. Hence, this interval is closed; let us denote it
$\Delta_{C^{(m)}} = [\delta_{m-}, \delta_{m+}]$. Then, the function $g(\delta)$ is
\begin{align}
\label{eq:gDelta}
g(\delta) = \rho^{(m)} + \delta (C_\text{req} -  C^{(m)}),
\; \delta \in [\delta_{m-}, \delta_{m+}]
\end{align}
where $\rho^{(m)}$ is the unique $\rho^\star(\delta)$ for $ \delta \in [\delta_{m-}, \delta_{m+}]$.

Now let $C^{(1)} = \max \{ C \in \mathcal{C} : C \leq C_\text{req}\}$ and $C^{(2)} = \min \{ C \in \mathcal{C} : C > C_\text{req}\}$ with associated reliability $\rho^{(1)}$ and $\rho^{(2)}$. Note that $C^{(1)} \leq C_\text{req} < C^{(2)}$.
Their associated $\delta$ range is
$\Delta_{C^{(1)}}=[\delta_{1-},\delta_{1+}]$ and
$\Delta_{C^{(2)}}=[\delta_{2-},\delta_{2+}]$. Furthermore, we have
$\delta^\star \triangleq \delta_{2+}=\delta_{1-}$ because $h(\delta)$ is a
continuous function.
Since $C^\star(\delta)$ is a non-increasing function over $\delta$,
we have $C^\star(\delta) \leq C^{(1)} \leq C_\text{req}$ for $\delta \geq \delta_{1+}$
and $C^\star(\delta) \geq C^{(2)} > C_\text{req}$
for $\delta \leq \delta_{2-}$.
Thus, we have $C^\star(\delta) \leq C_\text{req}$ for $\delta \geq \delta_{1-}$  and $C^\star(\delta) > C_\text{req}$ for $\delta \leq \delta_{2+}$.
Furthermore, by Eq.~\eqref{eq:gDelta},
$g(\delta)$ is a decreasing function for $\delta \leq \delta_{2+}$
since $C_\text{req} < C^{\star}(\delta)$
and a non-decreasing function for $\delta \geq \delta_{1-}$ since $C_\text{req} \geq C^{\star}$, so
the minimum value of $g(\delta)$ is obtained for
$\delta=\delta^\star=\delta_{2+}=\delta_{1-}$. The optimal $\delta^\star$ can be found from
$\rho^{(2)} + \delta^\star (C_\text{req} - C^{(2)} )
= \rho^{(1)} + \delta^\star (C_\text{req} - C^{(1)} ),$
and the maximum reliability (\textit{i.e.}, the minimum $g(\delta)$) is
\begin{align*}
\rho^\star=\rho^{(1)} +
\dfrac{C_\text{req}-C^{(1)}}{C^{(2)}-C^{(1)}} (\rho^{(2)} - \rho^{(1)}).
\end{align*}

Suppose the optimal policies to obtain $(\rho^{(1)}, C^{(1)})$ and $(\rho^{(2)},
C^{(2)})$ are $\pi^{(1)}$ and $\pi^{(2)}$ respectively.
Note that these two optimal policies for weighted sum maximization are history-independent and deterministic by~\cite{Put:94}.
The policy $\pi^\star$ that randomizes between $\pi^{(1)}$ and $\pi^{(2)}$ with
probabilities
\begin{align*}
\theta^{(1)} = \dfrac{C^{(2)}-C_\text{req}}{C^{(2)}-C^{(1)}}
 \; \text{ and } \;
\theta^{(2)} = \dfrac{C_\text{req}-C^{(1)}}{C^{(2)}-C^{(1)}}
\end{align*}
achieves this maximum reliability. Thus, it is an optimal policy, which concludes the proof.
\hfill $\blacksquare$

\begin{lemma}\label{lem:mono}
For operators $g_{\rho}(X)$ and $h_{\rho}(X)$, matrices $X\in\mathbb{S}_{>0}^{n}$ and $Y\in\mathbb{S}_{>0}^{n}$ and scalars $\rho,\;\rho^{\prime}\in[0,1]$, the following facts are true.
\begin{align*}
\renewcommand{\arraystretch}{1.2}
\setlength{\arraycolsep}{2.0pt}
\begin{array}{rrclcrcl}
	\mathbf{a)} & X & \leq & Y & \Longrightarrow & g_{\rho}(X) & \leq & g_{\rho}(Y) \;, \\
	\mathbf{b)} & X & \leq & Y & \Longrightarrow & h_{\rho}(X) & \leq & h_{\rho}(Y) \;, \\
  	\mathbf{c)} & \rho & \leq & \rho^{\prime} & \Longrightarrow & g_{\rho}(X) & \geq & g_{\rho^{\prime}}(X) \;, \\
  	\mathbf{d)} & \rho & \leq & \rho^{\prime} & \Longrightarrow & h_{\rho}(X) & \geq & h_{\rho^{\prime}}(X) \;.
\end{array}
\end{align*}
\end{lemma}
{\noindent \textit{Proof:}}
a) $g_{\rho}(X)\leq g_{\rho}(Y)$ holds as $g_{\rho}(X)$ is affine in $X$. The proof can be found in~\cite[Lemma~1(c)]{SSF+:04}.

b) As $h_{\rho}$ is a special form of $g_{\rho}$ by setting $\Phi=I$ and $R_{w}=0$, we immediately obtain the condition $h_{\rho}(X)\leq h_{\rho}(Y)$.

c) The detail proof is given in~\cite[Lemma~1(d)]{SSF+:04}.

d) As $h_{\rho}$ is a special form of $g_{\rho}$ by setting $\Phi=I$ and $R_{w}=0$, we easily conclude $h_{\rho}(X)\geq h_{\rho^{\prime}}(X)$ from (c).
\hfill $\blacksquare$

\begin{lemma}\label{lem:trace}
Suppose $Y\in\mathbb{S}_{>0}^{n}$ and $Z\in\mathbb{S}_{>0}^{n}$. If $Y\leq Z$, then $\mbox{Tr}(XY)\leq\mbox{Tr}(XZ),~\forall X\in\mathbb{S}_{>0}^{n}$.
\end{lemma}
{\noindent \textit{Proof:}}
Since $Y\leq Z$, it follows that $X^{\frac{1}{2}}YX^{\frac{1}{2}}\leq X^{\frac{1}{2}}ZX^{\frac{1}{2}}$. Since trace is a monotone function on the positive definite matrices, we get $\mbox{Tr}(XY)\leq \mbox{Tr}(XZ)$. \hfill $\blacksquare$

\bibliographystyle{IEEEtran}
\bibliography{codesign}

\end{document}